\documentclass[11pt]{article}
\usepackage{framed}
\usepackage{amsmath}
\usepackage{amsthm}
\usepackage{amsfonts}
\usepackage{amssymb}
\usepackage[utf8]{inputenc}
\usepackage{color}
\usepackage[T1]{fontenc}
\usepackage[margin=2.5cm]{geometry}
\usepackage{stmaryrd}
\usepackage{graphicx}
\usepackage{dsfont}
\usepackage{mathrsfs}
\usepackage{combelow}
\usepackage[font=small,labelfont=bf]{caption}
\usepackage{tikz}
\usetikzlibrary{patterns}
\usepackage{float}

\usepackage{xcolor}
\usepackage[hidelinks]{hyperref}
\hypersetup{
 colorlinks,
 linkcolor={red!50!black},
 citecolor={blue!50!black},
 urlcolor={blue!80!black}
}

\theoremstyle{plain}

\renewcommand{\Re}{\operatorname{Re}}

\newcommand{\Ind}{\mathds{1}}

\newcommand{\tend}[2]{\underset{#1\to #2}{\longrightarrow} }

\newcommand{\ntriple}[1]{{\left\vert\kern-0.25ex\left\vert\kern-0.25ex\left\vert #1 
    \right\vert\kern-0.25ex\right\vert\kern-0.25ex\right\vert}}

\newcommand{\der}[2]{\frac{\dd #1}{\dd #2}}

\newcommand{\dd}{\mathrm{d}}
\newcommand{\Diff}{\mathrm{D}}

\newcommand{\supp}{\mathrm{supp}\hspace{1mm}}

\newcommand{\eps}{\varepsilon}

\newcommand{\R}{\mathbb{R}}
\newcommand{\C}{\mathbb{C}}
\newcommand{\N}{\mathbb{N}}

\newcommand{\ds}{\displaystyle}

\renewcommand{\sp}{\hspace{0.2cm}}

\newcommand{\tgm}{\tilde{\gamma}}
\renewcommand{\tilde}{\widetilde}

\theoremstyle{plain}

\newtheorem{theo}{Theorem}[section]
\newtheorem{lemme}[theo]{Lemma}
\newtheorem{prop}[theo]{Proposition}

\newtheorem{remark}[theo]{Remark}

\newcommand\blfootnote[1]{%
  \begingroup
  \renewcommand\thefootnote{}\footnote{#1}%
  \addtocounter{footnote}{-1}%
  \endgroup
}

\title{Construction of unstable concentrated solutions of the Euler and gSQG equations}
\author{Martin Donati\footnote{Univ. Grenoble Alpes, Institut Fourier, F-38000
Grenoble, France. Contact : martin.donati@univ-grenoble-alpes.fr.}}
\date{}

\begin{document}

\maketitle


\begin{abstract}
    In this paper we construct solutions to the Euler and gSQG equations that are concentrated near unstable stationary configurations of point-vortices. Those solutions are themselves unstable, in the sense that their localization radius grows from order $\eps$ to order $\eps^\beta$ (with $\beta < 1$) in a time of order $|\ln\eps|$. This proves
    in particular that the logarithmic lower-bound obtained in previous papers (in particular [P. Buttà and C. Marchioro, \emph{Long time evolution of concentrated Euler flows with planar symmetry}, SIAM J. Math. Anal., 50(1):735–760, 2018]) about vorticity localization in Euler and gSQG equations is optimal. In addition we construct unstable solutions of the Euler equations in bounded domains concentrated around a single unstable stationary point. To achieve this we construct a domain whose Robin's function has a saddle point.
\end{abstract}

\section{Introduction}

We are interested in this paper in different active scalar equations from fluid dynamics: the two-dimensional incompressible Euler equations, used to describe an inviscid and incompressible fluid; and the Surface Quasi Geostrophic (SQG) equations, used as a model of geophysical flows. We also consider the generalized Surface Quasi-Geostrophic equations (gSQG) that interpolates between the Euler equations and the SQG equations. Let $\omega : \R^2 \to \R$ be the active scalar, which we will refer to as the \emph{vorticity} as in the Euler equations, and $u: \R^2 \to \R^2$ be the velocity of the fluid. Then the Euler equations, the SQG equations and the gSQG equations in the plane can all be written in the form
\begin{equation*}
\begin{cases}
    \partial_t \omega(x,t) + u \cdot \nabla \omega(x,t) = 0,& \forall (x,t) \in \R^2\times (0,+\infty) \\
    u(x,t) = \nabla^\perp\Delta^{-s}(\omega)(x,t), & \forall (x,t) \in \R^2\times [0,+\infty) \\
    u(x,t) \to 0,  & \text{as } |x| \to +\infty, \sp \forall t \in [0,+\infty) \\
    \omega(x,0) = \omega_0(x), & \forall x \in \R^2,
\end{cases}
\end{equation*}
with $s \in [\frac{1}{2},1]$. The Euler equations correspond to the case $s=1$, the SQG equations to the case $s = \frac{1}{2}$, and the gSQG equations are the family $s \in (\frac{1}{2},1)$. For details on the geophysical model see for instance \cite{Pedlowsky_1987,Vallis_2006}. Let us recall that the fundamental solution of $\Delta^s$ in the plane is given for $x \neq y$ by
\begin{equation*}
    G_s(x,y) = \begin{cases}
    \frac{1}{2\pi} \ln|x-y| & \text{ if $s = 1$}, \\
    \frac{\Gamma(1-s)}{2^{2s}\,\pi\,\Gamma(s)} \frac{1}{|x-y|^{2-2s}}  & \text{ if $\frac{1}{2} \le s < 1$},
    \end{cases}
\end{equation*}
where $\Gamma$ is the standard Gamma function. Therefore for every $s \in (\frac{1}{2}, 1]$ there exists a constant $C_s$ such that
\begin{equation*}
    \nabla_x G_s(x,y) = C_s \frac{x-y}{|x-y|^{4-2s}}.
\end{equation*}
This motivates us to define $\alpha := 3-2s$, so that $|\nabla_x G_\alpha(x,y)| \le C_\alpha |x-y|^{-\alpha}$, for $\alpha \in [1,2)$ and the appropriate choice of $C_{\alpha}$. In conclusion, the equations that we consider are the family of equations
\begin{equation}\label{eq:gSQG}
    \begin{cases}
    \partial_t \omega(x,t) + u(x,t)\cdot \nabla \omega(x,t) = 0, \vspace{2mm} \\
    \ds u(x,t) = \int_{\R^2} C_\alpha \frac{(x-y)^\perp}{|x-y|^{\alpha+1}} \omega(y,t)\dd y, \\
    \omega(x,0) = \omega_0(x).
\end{cases}
\end{equation}
We recall that the Euler case corresponds to $s = \alpha = 1$. We observe that we necessarily have that $\nabla \cdot u = 0$ when it makes sense. When $\alpha \ge 2$, the kernel cease to be $L^1_{\mathrm{loc}}$ so the Biot-Savart law in \eqref{eq:gSQG} does not make sense anymore. In that case the Biot-Savart law should be expressed differently, see for instance \cite{Marchand_2008}. In this paper we only consider $\alpha \in [1,2)$.

For the Euler equations, we have global existence and uniqueness of both strong solutions, and weak solutions in $L^1 \cap L^\infty$ from the Yudovitch theorem~\cite{YUDOVICH19631407}. 
When $\alpha \in (1,2)$, from \cite{Geldhauser_Romito_2020}, we know the existence, but not uniqueness, of weak solutions in $L^1 \cap L^\infty$ of equations~\eqref{eq:gSQG}. The existence of strong solutions is only known locally in time, and the blow-up of the SQG equations is an important open problem. 

\medskip

In this paper we construct families of particular initial data $\omega_0$ such that \emph{any} (since it may not be unique) solution of~\eqref{eq:gSQG} satisfies various constraints. These functions $\omega_0$ always lie in $L^1 \cap L^\infty$, but can also be taken $C^\infty$. Hence when $\alpha =1$, the solution remains $C^\infty$ for all time, but in general, we can only consider global in time $L^1 \cap L^\infty$ solutions of \eqref{eq:gSQG}. 

\vspace{2mm}

We now focus on the particular situation where the active scalar is \emph{concentrated} into small blobs as follows. We denote by $D(z,r)$ the disk of radius $r$ centered in $z$. We consider solutions $\omega$ satisfying
\begin{equation*}
    \omega = \sum_{i=1}^N \omega_{i} \sp \text{ and } \sp \supp \omega_{i} \subset D(z_i(t),r(t)),
\end{equation*}
with $r(t)$ being small in some sense. A classical model to describe concentrated solutions of equations~\eqref{eq:gSQG} is the point-vortex model. The principle is to approximate the blob $\omega_i$ by the Dirac mass $a_i\delta_{z_i(t)}$, where the intensity
\begin{equation*}
    a_i = \int_{\R^2} \omega_i(x,t)\dd x,
\end{equation*}
is constant in time\footnote{See \cite{Geldhauser_Romito_2020} Corollary 2.7 for the case of weak solutions. For strong solutions this is a direct consequence of the fact that $\nabla\cdot u =0$, see for instance \cite{Marchioro_Pulvirenti_1993}.}. The dynamics of those Dirac masses, that we call \emph{point-vortices}, is then given by the system
\begin{equation}\label{eq:PVSalpha}\tag{$\alpha$-PVS}
    \forall i \in \llbracket 1,N\rrbracket, \quad \der{}{t} z_i(t) = C_\alpha\sum_{\substack{1 \le j \le N \\ j\neq i}}  a_j \frac{(z_i(t)-z_j(t))^\perp}{|z_i(t)-z_j(t)|^{\alpha +1}}.
\end{equation}
This system of equations is often called $\alpha$-point-vortex system, or $\alpha$-model. We recall that this model is mathematically justified, see for instance \cite{ marchioro1993VorticiesAndLocalization,Turkington1987,Smets_VanSchaftingen_2010,marchiorobutta2018,Davila_DelPino_Musso_Wei_2020_Gluing} for the Euler case and \cite{Geldhauser_Romito_2020, Rosenzweig_2020, Godard-Cadillac_Gravejat_Smets_2020, Cavallaro_Garra_Marchioro_2021} for the gSQG case. In particular, on a finite time interval $[0,T]$, if no collisions of point-vortices occurs, then if 
\begin{equation*}
    \omega_{0,\eps} \tend{\eps}{0} \sum_{i=1}^N a_i\delta_{z_i(0)}
\end{equation*} 
weakly in the sense of measures, then 
\begin{equation*}
    \omega_\eps(t) \tend{\eps}{0} \sum_{i=1}^N a_i\delta_{z_i(t)},
\end{equation*} where the $t\mapsto z_i(t)$ are the solutions of the point-vortex dynamics. This means that point-vortices are a singular limit of solutions of the associated PDE. Fore a more detailed introduction of the point-vortex system, we refer the reader to \cite{Marchioro_Pulvirenti_1993}.

\section{Vorticity confinement and main results}

We now introduce the long time vorticity confinement problem, recall some important theorems on the subject and state our main results.

\subsection{Long time confinement problem}

We define the \emph{long time confinement problem} as the following -- see \cite{marchiorobutta2018}. Let $N \in \N^*$ and $\eps > 0$. For each $i \in \{1,\ldots,N\}$ let $z_i^* \in \R^2$ chosen pairwise distinct and $a_i \in \R^* = \R \setminus \{ 0\}$. Assume that $\omega_0$ is such that
\begin{equation}\label{hyp:omega0}
    \begin{cases}
    \ds \omega_0 = \sum_{i=1}^N \omega_{0,i} \sp \text{ and } \sp \supp \omega_{0,i} \subset D(z_i^*,\eps), \vspace{1mm}\\
    \ds \omega_{0,i} \text{ has a sign} \sp \text{ and } \sp \int_{\R^2} \omega_{0,i}(x) \dd x = a_i, \vspace{1mm}\\
    \ds |\omega_0| \le C\eps^{-\nu}, \quad \text{ for some } \nu \ge 2, \vspace{1mm}\\
    \ds |z_i^*(t) - z_j^*(t)| > 0, \quad \forall t \in [0,+\infty),
    \end{cases}
\end{equation}
were we denote by $z_i^*(t)$ the associated solution of the point-vortex dynamics with initial data $z_i^*(0) = z_i^*$ and intensities $a_i$. The last hypothesis ensures that the point-vortex dynamics has a global in time solution: no collision occurs. This is not a very restrictive hypothesis since it is known that the point-vortex system has a global solution for almost any initial data, in the sense of the Lebesgue measure. This was proved for the Euler point-vortex dynamics (namely equations \eqref{eq:PVSalpha} for $\alpha = 1$) in the torus \cite{DurrPulvirenti_1982}, in bounded domains \cite{Donati_2021} and in the plane\footnote{With an additional hypothesis on the intensities: $\sum_{i \in P} a_i \neq 0$, for any $P \subset \{1,\ldots,N\}$, with $P \neq \emptyset$. This hypothesis was then weakened in \cite{Godard-Cadillac_2022_Vortex_Collapses} to $P \neq \{1,\ldots,N\}$.} \cite{Marchioro_Pulvirenti_1984}. For the general $\alpha$-model \eqref{eq:PVSalpha} it was proved in the plane\footnote{With the same additional hypothesis.} \cite{Cavallaro_Garra_Marchioro_2013_Localization_active_scalar, Godard-Cadillac_2022_Vortex_Collapses}.

Let $\beta < 1$. We introduce the exit time:
\begin{equation*}
    \tau_{\eps,\beta} = \sup\left\{t\ge 0 \text{ such that } \forall s \in [0,t], \sp  \supp \omega(\cdot,s) \subset \bigcup_{i=1}^N D(z_i^*(s),\eps^\beta) \right\}.
\end{equation*} 
The long time confinement problem consists in obtaining a lower-bound on $\tau_{\eps,\beta}$ in order to describe how long the approximation of a concentrated solution of equations~\eqref{eq:gSQG} by the point-vortex model \eqref{eq:PVSalpha} remains valid. Results have been obtained in \cite{marchiorobutta2018,Donati_Iftimie_2020,Cavallaro_Garra_Marchioro_2021}. In the following, we recall some of them and state our main results, starting with the case $\alpha =1$.

\subsection{Result for Euler equations in the plane}
A first general result was obtained in \cite{marchiorobutta2018}.
\begin{theo}[Marchioro-Buttà, \cite{marchiorobutta2018}]\label{theo:MB}
Let $\beta < 1/2$. Then there exists $\xi_0 > 0$ such that for every $\eps>0$ small enough, for any $\omega_0$ satisfying~\eqref{hyp:omega0} for some $\nu \ge 2$, the solution $\omega$ of the Euler equations satisfies
 \begin{equation*}
    \tau_{\eps,\beta} > \xi_0 |\ln \eps|.
\end{equation*}
\end{theo}
In special cases, this can be improved. For instance, with the same hypotheses than those of Theorem~\ref{theo:MB}, but assuming furthermore that $N=1$, one can easily obtain that $\tau_{\eps,\beta} \ge \eps^{-\xi_0}$, for some $\xi_0 > 0$. An extension of this result is the following.

\begin{theo}[Marchioro-Buttà, \cite{marchiorobutta2018}]
Let $\beta < 1/2$. Then there exists $\xi_0> 0$ and a configuration of point-vortices, namely a choice of $a_i$ and $z_i^*$, such that for any $\eps > 0$ small enough and for any $\omega_0$ satisfying~\eqref{hyp:omega0}, the solution $\omega$ of the Euler equations satisfies
 \begin{equation*}
    \tau_{\eps,\beta} > \eps^{-\xi_0}.
\end{equation*}
\end{theo}
The configuration is a self-similar expanding configuration of three point-vortices. The key point is that as point-vortices move far from each other, their mutual influence decreases with time.

In this paper, we want to prove that the logarithmic bound obtained in Theorems~\ref{theo:MB} is optimal, namely that there are solutions of \eqref{eq:gSQG} satisfying \eqref{hyp:omega0} such that $\tau_{\eps,\beta} \le \xi_1 |\ln \eps|$. We prove the following result.
\begin{theo}\label{theo:N=3}
There exists $\beta_0 < 1/2$, $\nu \ge 2$ and a configuration $\big((z_i^*)_i,(a_i)_i\big)$ of point-vortices with $N=3$ such that for every $\beta \in (\beta_0,1)$, for any $\xi_1 > \frac{4\pi}{3}(1-\beta)$ and for every $\eps>0$ small enough there exists $\omega_0$ satisfying \eqref{hyp:omega0} such that the solution $\omega$ of the Euler equations in the plane satisfies
 \begin{equation*}
      \tau_{\eps,\beta} \le \xi_1 |\ln\eps|.
 \end{equation*}
\end{theo}
This confirms that the logarithmic bound obtained in Theorems~\ref{theo:MB} is optimal.

\subsection{Result for the gSQG equations in the plane}
A result similar to Theorem~\ref{theo:MB} has been obtained for the gSQG equations.
\begin{theo}[Cavallaro-Garra-Marchioro, \cite{Cavallaro_Garra_Marchioro_2021}]\label{theo:CGM}
Let $\alpha \in (1,2)$ and $\beta$ such that $0<\beta<\frac{4-2\alpha}{5-\alpha}<\frac{1}{2}$. Then there exists $\xi_0 > 0$ such that for every $\eps>0$ small enough, for any $\omega_0$ satisfying \eqref{hyp:omega0} with $\nu = 2$ and for any $\omega$ a weak solution of \eqref{eq:gSQG} we have that
\begin{equation*}
    \tau_{\eps,\beta} > \xi_0 |\ln \eps|.
\end{equation*}
\end{theo}
\begin{remark}
     Please note that in the hypotheses of Theorem~\ref{theo:CGM} is assumed that $\nu = 2$, which is not in the hypotheses of Theorem~\ref{theo:MB}. Actually, we claim that in their proof of Theorem~\ref{theo:CGM}, the authors of \cite{Cavallaro_Garra_Marchioro_2021} only need to assume the existence of $\nu \ge2$, and not $\nu = 2$. This is due to Lemma 2.4 of \cite{Cavallaro_Garra_Marchioro_2021}.
\end{remark}

We then prove the following.
\begin{theo}\label{theo:N=3sqg}
Let $\alpha \in [1,2)$. Then there exists $\nu \ge 2$, an initial configuration $\big((z_i^*)_i,(a_i)_i\big)$ of point-vortices with $N=3$, and $\beta_0 < \frac{4-2\alpha}{5-\alpha}$ such that for every $\beta \in (\beta_0,1)$, for every $\xi_1 > \frac{1-\beta}{C_\alpha(2-2^{-\alpha} )\sqrt{\alpha}}$  and for every $\eps>0$ small enough there exists $\omega_0$ satisfying \eqref{hyp:omega0} such that any solution $\omega$ of \eqref{eq:gSQG} satisfies
 \begin{equation*}
      \tau_{\eps,\beta} \le \xi_1 |\ln\eps|.
 \end{equation*}
\end{theo}
This confirms that the logarithmic bound obtained in Theorems~\ref{theo:CGM} is optimal. 

\begin{remark}
    Both in Theorems~\ref{theo:N=3} and~\ref{theo:N=3sqg}, the lower-bound for $\xi_1$ is not optimal. Moreover, $\omega_0$ is localized initially in a disk of size $\eps$ but the size of its support is of order $\eps^{\nu/2}$. In Appendix~\ref{sec:compConsts}, we give details how concentrated the initial data needs to be depending on the construction, and give examples constructions involving more point-vortices, which improve the bounds for $\xi_1$ and $\nu$.
\end{remark}

\subsection{Results for the Euler equations in bounded domains}
In a second part of this paper, we turn to a new situation. We focus on the Euler equations, namely the case $\alpha=1$, but in a bounded domain $\Omega$.  Let us recall the Euler equations in a bounded and simply connected domain $\Omega \subset \R^2$:
    \begin{equation}\label{eq:Eu}\tag{Eu}
    \begin{cases}
    \partial_t \omega + u \cdot \nabla \omega = 0, \vspace{2mm} \\
    \ds u =\nabla^\perp \Delta^{-1}\omega, \\
    u \cdot n = 0, & \text{ on } \partial\Omega, \\
    \omega(x,0) = \omega_0(x), & \text{ on } \Omega.
    \end{cases}
\end{equation}
When being far from the boundary, one can express the effect of the boundary as a Lipschitz exterior field. This trick makes it very easy to extend Theorem~\ref{theo:MB} to the case of bounded domains -- as it is suggested in \cite{marchiorobutta2018}.

In that same article, the authors proved that when the initial vorticity is concentrated near the center of a disk, namely that $\Omega= D(0,1)$, $N=1$ and $z_1 = 0$, then we obtain the same power-law lower-bound $\tau_{\eps,\beta} \ge \eps^{-\xi_0}$ than with expanding self-similar configurations. This result has been generalized to other bounded domains in \cite{Donati_Iftimie_2020}. This is due to a strong stability property induced by the shape of the boundary. Here we are interested in the opposite situation: we construct a domain whose boundary creates an instability.
We then obtain a third and final result, different from Theorems~\ref{theo:N=3} and \ref{theo:N=3sqg} because it only involves a single blob.
\begin{theo}\label{theo:boundedDomain}
There exists a smooth bounded domain $\Omega$ and $\beta_0 < 1/2$ such that for every $\beta \in (\beta_0,1)$, for every $\xi_1>(1-\beta)\frac{2\pi}{\sqrt{3}}$ and every $\eps >0$ small enough, there exists $\omega_0$ satisfying \eqref{hyp:omega0} with $N=1$ and $\nu = 4$ such that the solution $\omega$ of \eqref{eq:Eu} satisfies
 \begin{equation*}
      \tau_{\eps,\beta} \le \xi_1 |\ln\eps|.
 \end{equation*}
\end{theo}

We prove Theorems~\ref{theo:N=3}, \ref{theo:N=3sqg} and \ref{theo:boundedDomain} using the same plan: we construct a solution of the point-vortex dynamics that move away from its initial position exponentially fast, then construct a solution concentrated around it in the sense of hypothesis \eqref{hyp:omega0}.

The paper is organized as follows. In Section~\ref{sec:prelim} we give several definitions and expose in details the plan of the proofs and the main tools. In Section~\ref{sec:proof_SQG} we do the explicit construction to prove Theorems~\ref{theo:N=3} and \ref{theo:N=3sqg}. Finally in Section~\ref{sec:proofbounded} we prove Theorem~\ref{theo:boundedDomain}.

\section{Outline of the proofs}\label{sec:prelim}
In this section we expose the main tools required for the proofs of our results. Before going any further, let us introduce some notation.

In the rest of the paper,
\begin{itemize}
    \item $|z|$, for $z \in \R^2$ designates the usual $2$-norm, or modulus,
    \item $z^\perp = (-z_2,z_1)$,
    \item $|Z|_\infty$ for $Z =(z_1,\ldots,z_N) \in (\R^2)^N$ designates $\max_{1\le i \le N} |z_i|$,
    \item $C$ is a name reserved for constants whose value is not relevant, and may change from line to line,
    \item $D(z,r)$ is the disk (in $\R^2$) or radius $r$ centered in $z$.
\end{itemize}
Please notice that our construction -- in particular $\omega_0$ -- depends on $\eps$, though we do not write the dependence of each quantity in $\eps$ for the sake of legibility.

\subsection{Plan}

The proofs of Theorems \ref{theo:N=3}, \ref{theo:N=3sqg} and \ref{theo:boundedDomain} rely on two main steps. We first look for an \emph{unstable} stationary configuration of point-vortices. Then we control the behaviour of a well prepared solution initially concentrated around this configuration.

Let us give some details on each step.

\paragraph{Step 1: constructing an unstable vortex configuration.} 
\text{ }
\vspace{1mm}

We consider the dynamical system $\der{}{t} Z(t) = f(Z(t))$. Then we say that $Z^*$ is a stationary point of the dynamics if $f(Z^*) = 0$, and that it is \emph{unstable} if $\Diff f(Z^*)$ has an eigenvalue with positive real part.

At this step, we first aim to choose $N$, $\Omega$ (when necessary), the family of $a_i$ and $z_i^*$. The trick is following: we choose intensities $a_i \in \R^*$ and a point $Z^* =(z_1^*,\ldots,z_N^*) \in (\R^2)^N$  which is a stationary and unstable initial datum of the point-vortex dynamics~\eqref{eq:PVSalpha}. Let us notice that choosing a stationary configuration $Z^*$ ensures that the hypothesis $|z_i^*(t) - z_j^*(t)| > 0$ for all $t \ge 0$ is always satisfied. The consequence of the instability is that for every $\eps >0$, there exists an initial configuration of point-vortices $Z_0$ such that
\begin{equation*}
    |Z^* - Z_0|_\infty = \frac{\eps}{2},
\end{equation*}
and the solution $t\mapsto Z(t)$ such that $Z(0) = Z_0$ move away exponentially fast from $Z^*$, therefore implying that for any $\beta < 1$,
\begin{equation*}
    \tau_Z := \sup\left\{t\ge 0 \text{ such that } \forall s \in [0,t] \; , \; |Z(s)-Z^*|_\infty \le 2\eps^\beta \right\} \le \xi_0 |\ln\eps|.
\end{equation*}
This problem is much simpler than the original one since we are investigating the behaviour of solutions of a system of ordinary differential equations, the point-vortex dynamics, instead of a solution of a partial derivative equation. More precisely, we have the following proposition, obtained as a corollary of Theorem 6.1, Chapter 9 of \cite{hartman1982ordinary}, and proved in Section~\ref{sec:proof_prelim}.

\begin{prop}\label{prop:solution_divergente}
Let $f : (\R^2)^N \to (\R^2)^N$. We consider the differential equation 
\begin{equation}\label{eq:diff_eq_generale}
    \der{}{t} Z(t) = f(Z(t)).
\end{equation}
Assume that there exists $Z^* \in (\R^2)^N$ such that $f(Z^*) = 0$. Assume furthermore that $\Diff f(Z^*)$ has an eigenvalue with positive real part $\lambda_0 > 0$. 

Then for any $\lambda < \lambda_0$, for every $\eps >0$ small enough and for any $\beta \in (0,1)$ there exists and a choice of $Z_0$ such that $|Z_0 - Z^*|_\infty  = \eps/2$, $$\tau_Z \le \frac{1-\beta}{\lambda} |\ln \eps|.$$ 
\end{prop}

In conclusion, proving that $\tau_Z \le \xi_1 |\ln \eps|$ simply relies on finding an eigenvalue with positive real part of the Jacobian matrix of the dynamic's functional.

\paragraph{Step 2: constructing the approximation} 
\text{ }
\vspace{1mm}

The idea is then to prove that a solution $\omega$ with well prepared initial data $\omega_0$ satisfying \eqref{hyp:omega0} satisfies that for every $t \le \tau_{\eps,\beta}$,
\begin{equation*}
    |B(t)-Z(t)|_\infty = o(\eps^\beta),
\end{equation*}
where
\begin{equation*}
    B_i(t) =  \frac{1}{a_i}\int_{\R^2} x \omega_i(x,t) \dd x,
\end{equation*}
and $B(t) = (B_1(t),\ldots,B_N(t))$. The conclusion then comes from the fact that by construction (Step 1), there exists $t_1 \le \frac{1-\beta}{\lambda}$ such that $|Z(t_1)-Z^*|_\infty = 2\eps^\beta$ and thus for $\eps$ small enough, $\tau_{\eps,\beta} \le t_1 \le \xi_1 |\ln \eps|$. 

In order to obtain this control on $|B(t)-Z(t)|_\infty$, we need to estimate the moment of inertia
\begin{equation*}
    I_i(t) = \frac{1}{a_i}\int_{\R^2} |x-B_i(t)|^2\omega_i(x,t) \dd x
\end{equation*}
of each blob. The constant $\nu$ in \eqref{hyp:omega0} intervenes when estimating $I(0)$. The critical part in this step is the competition between the growth of $I_i$, which loosen the control on $|B(t)-Z(t)|_\infty$, with the growth of $|Z(t)-Z^*|_\infty$. When we do not have constraints on $\nu$, then one can always choose $\nu$ large enough so that $I_i$ remains small long enough. However the difficulty arises when wanting to construct $\omega_0$ with $\nu =4$. This requires to be able to estimate precisely the growth of each $I_i$, of $|Z(t)-Z^*|_\infty$ and of $|B(t)-Z(t)|_\infty$.

All of this is captured in Theorem~\ref{prop:argFinal} presented in Section~\ref{sec:arg_final}.

\subsection{Confinement around an unstable configuration}\label{sec:arg_final}
Once the unstable configuration of point-vortices is obtained in Step 1, most of the work needed in the second step does not depend on that configuration nor on the specific framework. Therefore, we establish a general theorem that we will be able to apply for any suitable configuration of vortices, also including when appropriate the presence of a boundary.

To understand better the dynamics of each blob, we describe the influence of the other blobs or of the boundary by an an exterior field.
We assume that each blob is a solution of a problem
 \begin{equation}\label{eq:systemDuPB}
\begin{cases}
    \partial_t \omega_i(x,t) + \big(u_i(x,t) + F_i(x,t)\big)\cdot \nabla \omega_i(x,t) = 0, \vspace{2mm} \\
    \ds u_i(x,t) = \int_{\R^2} C_\alpha \frac{(x-y)^\perp}{|x-y|^{\alpha+1}}\omega_i(y,t)\dd y, \vspace{1mm}\\
    \omega_i(x,0) = \omega_{i,0}(x),\\
\end{cases}
\end{equation}
where $F_i$ is an exterior field that satisfies $\nabla \cdot F_i = 0$. 
Let $Z^* \in (\R^2)^N$ and $f \in C^1\Big( (\R^2)^N \, , \,  (\R^2)^N\Big)$ such that $f(Z^*)=0$. We write $f = (f_1,\ldots,f_N)$. For any $Z_0$, let $t\mapsto Z(t)$ be the solution of the problem
\begin{equation}\label{eq:pbz}
\begin{cases}
    \ds \der{}{t} Z(t) = f(Z(t)) \vspace{1mm}\\
    Z(0) = Z_0.
\end{cases}
\end{equation}
In this particular setting, for any $\beta \in (0,1)$, we have
\begin{equation*}
    \tau_{\eps,\beta} = \sup\left\{t\ge 0  \text{ such that } \forall s \in [0,t], \sp  \supp \omega(\cdot,s) \subset \bigcup_{i=1}^N D(z_i^*,\eps^\beta) \right\}.
\end{equation*}
Assuming that $|Z_0 - Z^*|_\infty= \eps/2$, we recall that
\begin{equation*}
   \tau_Z :=  \sup\left\{t\ge 0  \text{ such that } \forall s \in [0,t] \; , |Z(s) - Z^*|_\infty < 2\eps^\beta \right\}.
\end{equation*}

We then have the following theorem.
\begin{theo}\label{prop:argFinal}
Let $N \in \N^*$, $a_i \in \R^*$ for every $i \in \{1,\ldots,N\}$, $Z^* \in (\R^2)^N$. Let $f \in C^1\big((\R^2)^N,(\R^2)^N\big)$ and $F_i$ such that $\nabla\cdot F_i = 0$.
We assume the following.
\begin{itemize}
    \item [$(i)$] $f(Z^*) = 0$ and $\Diff f(Z^*)$ has an eigenvalue with positive real part $\lambda_0$,
    \item[$(ii)$] There exists $C$ such that for all $i\in\{1,\ldots,N\}$, $\forall t \le \tau_{\eps,\beta}$,
    \begin{equation*}
    \left|F_i(B_i(t),t) - f_i(B(t))\right| \le C \sum_{j=1}^N\sqrt{I_j},
\end{equation*}
    \item[$(iii)$] There exists constants $\kappa_0$, $\kappa_1$ and $\kappa_2$ such that $\forall i \in \{1,\ldots,N\}$, $\forall\, x,x' \in D(z_i^*,\eps^\beta)$,  $\forall\, t \le \tau_{\eps,\beta}$,
    \begin{equation}\label{eq:Flip}
        |F_i(x,t)-F_i(x',t)| \le \kappa_0|x-x'|,
    \end{equation}
    and
\begin{equation}\label{eq:Lip_F(kappa1)}
     \Big|(x-x')\cdot \big(F_i(x,t)-F_i(x',t)\big)\Big| \le \kappa_1 |x-x'|^2,
\end{equation}
and $\forall\, X,X' \in (\R^2)^N$ such that $|X-Z^*|_\infty \le 2\eps^\beta$ and $|X'-Z^*|_\infty\le 2\eps^\beta$,
\begin{equation}\label{eq:Lip_f(kappa2)}
      |f(X)-f(X')|_\infty \le \kappa_2 |X-X'|_\infty.
\end{equation}
\end{itemize}
Then there exists $\nu\ge 2$ and $\beta_0<\frac{4-2\alpha}{5-\alpha}$ such that for all $\beta\in (\beta_0,1)$, for every $\xi > \frac{1-\beta}{\lambda_0}$ and for every $\eps>0$ small enough, there exists $\omega_0$ satisfying \eqref{hyp:omega0} such that any $\omega = \sum_{i=1}^N \omega_i$ solution of the problem \eqref{eq:systemDuPB} for every $i$ satisfies
    $$\tau_{\eps,\beta} \le \xi |\ln \eps|.$$
\end{theo}


\begin{proof}

First, we use Hypothesis $(i)$ to apply Proposition~\ref{prop:solution_divergente} and get that for every $\eps$ small enough, there exists $Z_0$ such that $|Z_0-Z^*|_\infty = \eps/2$ and for every $\beta \in (0,1)$, the solution $Z$ of the problem~\eqref{eq:pbz} satisfies
\begin{equation}\label{eq:ttau}
    \tau_Z \le \frac{1-\beta}{\lambda} |\ln \eps|.
\end{equation}
Now let $\omega_0$ satisfying \eqref{hyp:omega0} for some $\nu \ge 2$ and such that
 \begin{equation}\label{eq:B(0)}
        B(0) = Z_0,
    \end{equation}
    and
    \begin{equation}\label{eq:I(0)}
        \forall i \in \{1,\ldots,N\}, \quad I_i(0) \le \eps^{\nu}.
    \end{equation}
This is always possible as stated in Remark~\ref{prop:choixOmega0} given in Appendix~\ref{sec:appendix_proofs}. Let $\omega = \sum_{i=1}^N \omega_i$ such that each $\omega_i$ is a solution of the problem~\eqref{eq:systemDuPB}.

We observe that if $\tau_{\eps,\beta} \le \tau_Z$, then we have the desired result. So for the sake of contradiction, we can assume that $\tau_{\eps,\beta} > \tau_Z$.

Recalling that $\omega_i$ solves \eqref{eq:systemDuPB}, we have that
\begin{equation*}
    \der{}{t} B_i = \frac{1}{a_i}\int F_i(x,t)\omega_i(x,t) \dd x
\end{equation*}
and
\begin{equation*}
    \der{}{t} I_i = \frac{2}{a_i}\int(x-B_i(t))\cdot(F_i(x,t)-F_i(B_i(t),t)) \omega_i(x,t)\dd x.
\end{equation*}
Indeed, if $\omega_i$ is smooth (when $\alpha = 1$ or before a possible regularity blow-up if $\alpha > 1$), these are classical computations. In general, $\omega_i \in L^1\cap L^\infty$ and these relations hold in the weak sense, see for instance \cite[Corollary 2.8]{Geldhauser_Romito_2020} or \cite{Cavallaro_Garra_Marchioro_2021}.

Now using Hypothesis $(iii)$, and observing from~\eqref{hyp:omega0} that $\frac{\omega_i}{a_i} \ge 0$, we get that
\begin{equation*}
     \left|\der{}{t} I_i\right| \le 2 \kappa_1 \int|x-B_i(t)|^2 \frac{\omega_i(x,t)}{a_i}\dd x = 2\kappa_1 I_i(t).
\end{equation*}
Therefore, we get that
\begin{equation}\label{eq:majI}
    I_i(t) \le I_i(0) e^{2\kappa_1 t}.
\end{equation}
We now want to estimate $|B(t)-Z(t)|_\infty$. For every $i \in \{1,\ldots,N\}$ we have that
\begin{align*}
    \left|\der{}{t}B_i(t)-f_i(B(t))\right| & = \left|\frac{1}{a_i}\int \big(F_i(x,t)-f_i(B(t))\big)\omega_i(x,t) \dd x \right|\\
    & = \left|\frac{1}{a_i}\int \big(F_i(x,t)-F_i(B_i(t),t) + F_i(B_i(t),t) - f_i(B(t))\big)\omega_i(x,t) \dd x \right|\\
    & \le \kappa_0 \int |x-B_i(t)|\frac{\omega_i(x,t)}{a_i} \dd x + C \sum_{j=1}^N \sqrt{I_j(t)},
\end{align*}
where we used hypotheses $(ii)$ and $(iii)$.
By the Cauchy Schwartz inequality we have that
\begin{equation*}
    \int |x-B(t)|\frac{\omega_i(x,t)}{a_i} \dd x = \left(\int |x-B(t)|^2\frac{\omega_i(x,t)}{a_i} \dd x\right)^{1/2} \left( \int \frac{\omega_i(x,t)}{a_i} \dd x \right)^{1/2}\le \sqrt{I_i(t)},
\end{equation*}
and therefore, since the result is now uniform in $i$,
\begin{equation*}
    \left|\der{}{t}B(t)-f(B(t))\right|_{\infty} \le C \sum_{j=1}^N\sqrt{I_j(t)}.
\end{equation*}
The value of $C$ is irrelevant and changes from line to line. The value of $\kappa_0$ is also irrelevant in the end and is absorbed in $C$.

We recall that by relation \eqref{eq:B(0)}, $B(0) = Z_0 = Z(0)$ and that $Z$ is a solution of the problem \eqref{eq:pbz} with $f$ being a Lipschitz map by Hypothesis $(iii)$. We now use a variant of the Gronwall's inequality --~Lemma~\ref{lemme:gronwall} given in appendix -- to obtain that
\begin{equation*}
    \big| B(t) - Z(t) \big|_\infty \le Ce^{\kappa_2 t} \int_0^t \sum_{j=1}^N\sqrt{I_j(s)}\dd s .
\end{equation*}
Using relation \eqref{eq:majI}, we have that
\begin{equation*}
    \sum_{j=1}^N\int_0^t \sqrt{I_j(s)}\dd s \le \sum_{j=1}^N\int_0^t \sqrt{I_j(0)}e^{\kappa_1 s}\dd s \le \sum_{j=1}^N\frac{\sqrt{I_j(0)}}{\kappa_1}e^{\kappa_1 t}. 
\end{equation*}
Recalling that by relation~\eqref{eq:I(0)}, $I_j(0) \le \eps^{\nu}$, we obtain that
\begin{equation*}
    \big| B(t) - Z(t) \big|_\infty \le C \eps^{\nu/2-(\kappa_1+\kappa_2) t|\ln \eps|}.
\end{equation*}
We are interested in proving that $\forall t \le \tau_Z$,
\begin{equation}\label{eq:close}
    \big| B(t) - Z(t) \big|_\infty = o(\eps^\beta),
\end{equation} 
as $\eps \to 0$.
For $\eps <1$ and for all $t \le \tau_Z$, by relation \eqref{eq:ttau}, for any $\lambda < \lambda_0$,
\begin{equation*}
    -t |\ln \eps| \ge \frac{1-\beta}{\lambda}.
\end{equation*}
Therefore, one sufficient condition to obtain \eqref{eq:close} is
\begin{equation}\label{eq:compcst}
    \nu > 2\frac{\kappa_1+\kappa_2}{\lambda_0}(1-\beta) + 2\beta,
\end{equation}
since in that case, one can choose $\lambda$ close enough to $\lambda_0$ such that
\begin{equation*}
    \nu/2 - \frac{\kappa_1+\kappa_2}{\lambda}(1-\beta) > \beta.
\end{equation*}
We observe that necessarily, $\kappa_2 \ge \lambda_0$ so $\frac{\kappa_1+\kappa_2}{\lambda_0} > 1$ and therefore the map $\phi : \beta \mapsto 2\frac{\kappa_1+\kappa_2}{\lambda_0}(1-\beta) + 2\beta$ is decreasing on $[0,1]$. Let 
\begin{equation}\label{eq:nu}
 \nu >  \phi\left(\frac{4-2\alpha}{5-\alpha}\right) =  \frac{2}{5-\alpha} \left( (1+\alpha)\frac{\kappa_1+\kappa_2}{\lambda_0} + 4-2\alpha\right).    
\end{equation} 
Since $\phi$ is also continuous, then there exists $\beta_0 < \frac{4-2\alpha}{5-\alpha}$ such that $\forall \beta \in (\beta_0,1)$, $$\nu > \phi(\beta).$$ In particular, for every $\beta \in (\beta_0,1)$, relation~\eqref{eq:compcst} holds true and thus relation \eqref{eq:close} too.

We conclude by observing that for $\eps$ small enough,
\begin{equation*}
    |B(\tau_Z) - Z^*|_\infty \ge |Z(\tau_Z) - Z^*|_\infty - |B(\tau_Z) - Z(\tau_Z)|_\infty \ge 2\eps^\beta - o(\eps^\beta) \ge \eps^\beta.
\end{equation*}
Since $\forall\, t < \tau_{\eps,\beta}$, we have that $|B(t)-Z^*|_\infty < \eps^\beta$ by definition of $\tau_{\eps,\beta}$, the previous relation is in contradiction with $\tau_{\eps,\beta} > \tau_Z$. Therefore,
\begin{equation*}
    \tau_{\eps,\beta} \le \tau_Z \le \frac{1-\beta}{\lambda} |\ln \eps|,
\end{equation*}
for any $\lambda<\lambda_0$. In particular, for every $\beta \in (\beta_0,1)$, for every $\xi_1 > \frac{1-\beta}{\lambda_0}$, we proved that for $\eps$ small enough,
\begin{equation*}
    \tau_{\eps,\beta} \le \xi_1 |\ln \eps|.
\end{equation*}

\end{proof}


Let us mention that Hypothesis $(i)$ is a consequence of the choice of $Z^*$ and the $a_i$ made at Step 1, and that Hypothesis $(ii)$ is a consequence of the nature of the relation between the exterior fields $F_i$ and the map $f$. Of course, the result could not be true if the $F_i$ and $f$ were not related to each other.

The existence of $\kappa_0$, $\kappa_1$ and $\kappa_2$ in the Hypothesis $(iii)$ is a direct consequence of the fact that every $F_i$ and $f$ are Lipschitz maps, which will always be the case in our framework. Even though relation~\eqref{eq:Lip_F(kappa1)} is a consequence of \eqref{eq:Flip}, only $\kappa_1$ and $\kappa_2$ intervene in the condition~\eqref{eq:nu} on $\nu$. This is important when one wants to compute precisely the bounds on $\nu$, as we do in Section~\ref{sec:proofbounded}, and for several other examples in Appendix~\ref{sec:compConsts}.

Finally, we notice that we actually proved that any solution $\omega$ starting from \emph{any} initial data $\omega_0$ satisfying \eqref{hyp:omega0} and relations~\eqref{eq:B(0)} and~\eqref{eq:I(0)} satisfies $\tau_{\eps,\beta} \le \xi_1|\ln \eps|$. Therefore, there is a large family of initial data, in terms for instance of the shape of the blob, that solves our problem.

\subsection{Proof of Proposition \ref{prop:solution_divergente}.}\label{sec:proof_prelim}

Let us recall Theorem 6.1, Chapter 9 of \cite{hartman1982ordinary}.
\begin{theo}\label{theo:center_manifold}
Let $f : (\R^2)^N \to (\R^2)^N$. We consider the differential equation 
\begin{equation*}
    \der{}{t} Z(t) = f(Z(t)).
\end{equation*}
Assume that there exists $Z^* \in (\R^2)^N$ is such that $f(Z^*) = 0$. Assume furthermore that $\Diff f(Z^*)$ has an eigenvalue with positive real part $\lambda_0 > 0$. Then there exists a solution of \eqref{eq:diff_eq_generale} such that $Z(t)$ exists some fixed neighbourhood of $Z^*$, that
\begin{equation*}
    Z(t) \tend{t}{-\infty} Z^*,
\end{equation*}
and that
\begin{equation*}
    \frac{1}{t}\ln |Z(t)-Z^*|_\infty \tend{t}{-\infty} \lambda_0.
\end{equation*}

\end{theo}

We now prove Proposition~\ref{prop:solution_divergente}. Let $\beta \in (0,1)$ and let $\tilde{Z}$ a solution of \eqref{eq:diff_eq_generale} given by Theorem \ref{theo:center_manifold}. Since $\tilde{Z} \tend{t}{-\infty} Z^*$ and since $\tilde{Z}$ exits some fixed neighbourhood of $Z^*$, for $\eps$ small enough, there exists $t_0$ and $t_1$ such that
\begin{equation*}
    \begin{cases}
    -\infty < t_0 < t_1 \\
    t_1 \to -\infty & \text{ as } \eps \to 0\\
    |\tilde{Z}(t_1)-Z^*|_\infty = 2\eps^\beta \\
    |\tilde{Z}(t_0)-Z^*|_\infty = \eps/2.
    \end{cases}
\end{equation*}
Let $Z(t) = \tilde{Z}(t + t_0)$, we have that $|Z(0) - Z^*|_\infty = \eps/2$ and that $\tau_Z \le t_1-t_0$ since $Z(t_1-t_0) = 2 \eps^\beta$. Moreover, since 
\begin{equation*}
    \ln |Z(t)-Z^*|_\infty = \lambda_0 t + o_{t \to -\infty}(t),
\end{equation*}
then for any $\eta \in (0,1)$, for $-t$ big enough we have that
\begin{equation*}
    1 - \eta < \frac{\ln |Z(t)-Z^*|_\infty}{\lambda_0 t} < 1 + \eta 
\end{equation*}
Therefore, for $\eps$ small enough, applying in $t_0$ and $t_1$ (we recall that $t_1 \to -\infty$ as $\eps \to 0$) we have that 
\begin{equation*}
    t_1 < \frac{\ln |Z(t_1)-Z^*|_\infty}{\lambda_0(1+\eta)} = \frac{\beta\ln \eps + \ln 2}{\lambda_0 (1+\eta)} = \frac{- \beta|\ln \eps|+\ln 2}{\lambda_0 (1+\eta)}
\end{equation*}
and
\begin{equation*}
    -t_0 < -\frac{\ln |Z(t_0)-Z^*|_\infty}{\lambda_0(1-\eta)} = \frac{|\ln\eps|+\ln 2}{\lambda_0 (1-\eta)}
\end{equation*}
and thus
\begin{equation*}
    t_1 - t_0 < |\ln\eps| \left( \frac{1}{\lambda_0(1-\eta)} - \frac{\beta}{\lambda_0(1+\eta)} + \frac{\ln 2 }{\lambda_0|\ln \eps| } \left( \frac{1}{1+\eta} + \frac{1}{1-\eta}\right) \right).
\end{equation*}
Therefore, by letting $\eta \to 0$, for any $\lambda < \lambda_0$, for $\eps$ small enough, 
\begin{equation*}
    t_1 - t_0 \le \frac{1-\beta}{\lambda} |\ln \eps|.
\end{equation*}
By definition, $\tau_Z \le t_1 - t_0 \le \frac{1-\beta}{\lambda} |\ln \eps|$.
This concludes the proof.

\section{Unstable configurations of multiple vortices  in the plane}\label{sec:proof_SQG}

In this section we show that we can choose a configuration $Z^*$ and intensities $a_i$ such that we can apply Theorem~\ref{prop:argFinal} to estimate the exit time $\tau_{\eps,\beta}$ of a solution of the equations~\eqref{eq:gSQG}. Please notice that Theorem~\ref{theo:N=3} is a direct consequence of Theorem~\ref{theo:N=3sqg} by taking $\alpha = 1$. Therefore, we work with some general $\alpha \in [1,2)$.

We start by introducing a family of point-vortex configurations $Z^*$ for any $N \ge 3$. We then construct the exterior fields $F_i$ such that the blobs $\omega_i$ is a solution of the problem \eqref{eq:systemDuPB}, and start proving each of the conditions to apply Theorem~\ref{prop:argFinal}. We then prove Theorem~\ref{theo:N=3}, and Proposition~\ref{prop:N=3} by computing explicitly the properties of our construction with $N=3$. In Appendix~\ref{sec:compConsts}, we use the construction with $N=7$ and $N=9$.

\subsection{Vortex crystals}\label{sec:vortex_crystals}

Let us introduce a family of stationary solutions of the $\alpha$-point-vortex model \eqref{eq:PVSalpha} that are part of the so called \emph{vortex crystals} family. For a more general study of vortex crystals and their stability, we refer the reader to \cite{AREF_2003_Vortex_Crystals,Morikawa_Swenson_1971_Interacting_Motion_Vortices,Cabral_Schmidt_2000_Stability_N+1_vortex}. Fore the sake of legibility, we identify $\C = \R^2$ for the position of point-vortices. We use the notation $x = \big( (x)_p,(x)_q \big) = (x)_p + i (x)_q$.

Let $N \ge 3$, and
\begin{equation*}
    K_\alpha(x,y) =  C_\alpha\frac{(x-y)^\perp}{|x-y|^{\alpha +1}},
\end{equation*}
so that the point-vortex equation \eqref{eq:PVSalpha} becomes
\begin{equation*}
    \forall i \in \llbracket 1,N\rrbracket, \quad \der{}{t} z_i(t) = \sum_{\substack{1 \le j \le N \\ j\neq i}}  a_j K_\alpha\big(z_i(t),z_j(t)\big).
\end{equation*}
In particular, by setting for all $Z=(z_1,\ldots,z_N)$:
\begin{equation*}
    f(Z) = \left( \sum_{j\neq i} a_j K_\alpha(z_i,z_j)\right)_{1\le i \le N},
\end{equation*}
then we have that
\begin{equation*}
    \der{}{t} Z(t) = f(Z(t)).
\end{equation*}

We consider $N$ point-vortices in the following configuration. The first $N-1$ points form a regular $(N-1)$-polygon, and the $N$-th vortex is placed at the center. For instance, by letting $\zeta = e^{i \frac{2\pi}{N-1}}$, where here $i$ denotes the complex unit, we set
\begin{align*}
    \forall j \in \{1,\ldots,N-1\}, \sp z_j^* = \zeta^{j}, & \quad   a_j = 1 \\
    z_N^* = 0, & \quad a_N \in \R.
\end{align*}
Then, (see for instance \cite{AREF_2003_Vortex_Crystals}) the solution of \eqref{eq:PVSalpha} with initial configuration $Z^* = (z_1^*,\ldots,z_N^*)$ satisfies $z_j(t) = e^{i \mu t}z_j^*$ for some angular velocity $\nu$ that does not depend on $j$. The motion of the whole configuration is a rigid rotation around 0, which makes it a so called vortex crystal.

This stands for any choice of $a_N \in \R$. Now we make a particular choice. For every $N \ge 3$, there exists $a_N \neq 0$ in the previous configuration such that the solution is stationary ($\mu = 0$). Indeed, let us compute the velocity of any point vortex (except the one at the center that is always stationary), for instance $z_{N-1} = \zeta^{N-1} = 1$.
\begin{align*}
    \der{}{t} z_{N-1}(0) & = \sum_{j=1}^{N-2} K_\alpha(1,\zeta^{j}) + a_N K_\alpha(1,0) \\
    & = C_\alpha\left( \sum_{j=1}^{N-2}\frac{(1-\zeta^j)}{|1-\zeta^j|^{\alpha +1}} + a_N \right)^\perp.
\end{align*}
Since $\overline{\zeta^j} = \zeta^{N-1-j}$, the quantity $\sum_{j=1}^{N-2}\frac{(1-\zeta^j)}{|1-\zeta^j|^{\alpha +1}}$ is a non vanishing real number and thus letting 
\begin{equation*}
    a_N = - \sum_{j=1}^{N-2}\frac{(1-\zeta^j)}{|1-\zeta^j|^{\alpha +1}}
\end{equation*}
enforces that $\der{}{t} z_{N-1}(0) = 0$. By symmetry, $\der{}{t} z_{j}(0) = 0$ for every $j \in {1,\ldots,N-1}$. As for $z_N$, it is always stationary, again by a symmetry argument, or by a simple computation.
\medskip

In conclusion, for any $N\ge 3$, we constructed a $N$-vortex configuration that is stationary.
In order to study the stability of the equilibrium $Z^*$, we compute $\Diff f(Z^*)$. To this end, let us compute
\begin{align*}
\frac{1}{C_\alpha}&  \big(K_\alpha(z_i + x, z_j + y) - K_\alpha(z_i,z_j) \big) \\ 
    & = \frac{\big((z_i+x)-(z_j+y)\big)^\perp}{\big|(z_i+x)-(z_j+y)\big|^{\alpha+1}} - \frac{(z_i-z_j)^\perp}{|z_i-z_j|^{\alpha+1}} \\
    & =  \frac{(z_i-z_j)^\perp}{\big|(z_i+x)-(z_j+y)\big|^{\alpha+1}} - \frac{(z_i-z_j)^\perp}{|z_i-z_j|^{\alpha+1}} + \frac{(x-y)^\perp}{\big|(z_i+x)-(z_j+y)\big|^{\alpha+1}} \\
    & = \frac{(z_i-z_j)^\perp}{|z_i-z_j|^{\alpha+1}} \left( 1-(\alpha+1)(x-y)\cdot\frac{z_i-z_j}{|z_i-z_j|^2}+o\big(|x|+|y|\big)-1\right) +\frac{(x-y)^\perp}{|z_i-z_j|^{\alpha+1}} + o\big(|x|+|y|\big),
\end{align*}
and finally, we obtain that
\begin{multline}\label{eq:developpmentK}
    \frac{1}{C_\alpha}  \big(K_\alpha(z_i + x, z_j + y) - K_\alpha(z_i,z_j) \big) \\ = -\frac{(z_i-z_j)^\perp}{|z_i-z_j|^{\alpha+3}}(\alpha+1)(x-y)\cdot(z_i-z_j)+\frac{(x-y)^\perp}{|z_i-z_j|^{\alpha+1}} + o\big(|x|+|y|\big).
\end{multline}
Since each coordinate $z_i$ of $f$ is of dimension two, we need some clarification on the notations. We now think of $f$ as the map 
\begin{equation*}
    \tilde{f} : \begin{cases}\R^{2N} \to \R^{2N} \\ \big(p_1,q_1,\ldots,p_N,q_N\big) \mapsto \big(\tilde{f}_{p_1},\tilde{f}_{q_1},\ldots,\tilde{f}_{p_N},\tilde{f}_{q_N}\big) =  f\big(p_1+iq_1,\ldots,p_N+iq_N\big).\end{cases}
\end{equation*}
Relation \eqref{eq:developpmentK} yields for $i \neq j$ that
\begin{equation*}
    \begin{cases}
        \ds \frac{1}{a_j C_\alpha}\partial_{p_j} \tilde{f}_{p_i} = -\frac{(z_i-z_j)_q}{|z_i-z_j|^{\alpha+3}}(\alpha+1)(z_i-z_j)_p \vspace{1mm}\\
        \ds \frac{1}{a_j C_\alpha}\partial_{p_j} \tilde{f}_{q_i} = \frac{(z_i-z_j)_p}{|z_i-z_j|^{\alpha+3}}(\alpha+1)(z_i-z_j)_p-\frac{1}{|z_i-z_j|^{\alpha+1}} \vspace{1mm}\\
        \ds \frac{1}{a_j C_\alpha}\partial_{q_j} \tilde{f}_{p_i} = -\frac{(z_i-z_j)_q}{|z_i-z_j|^{\alpha+3}}(\alpha+1)(z_i-z_j)_q +   \frac{1}{|z_i-z_j|^{\alpha+1}} \vspace{1mm}\\
        \ds \frac{1}{a_j C_\alpha}\partial_{q_j} \tilde{f}_{q_i} = \frac{(z_i-z_j)_p}{|z_i-z_j|^{\alpha+3}}(\alpha+1)(z_i-z_j)_q,
    \end{cases}
\end{equation*}
and for $i = j$,
\begin{equation*}
    \begin{cases}
        \ds \partial_{p_i} \tilde{f}_{p_i} =  -\sum_{j\neq i} \partial_{p_j} \tilde{f}_{p_i}\vspace{1mm}\\
        \ds \partial_{p_i} \tilde{f}_{q_i} =  -\sum_{j\neq i} \partial_{p_j} \tilde{f}_{q_i} \vspace{1mm}\\
        \ds \partial_{q_i} \tilde{f}_{p_i} =  -\sum_{j\neq i}\partial_{q_j} \tilde{f}_{p_i}\vspace{1mm}\\
        \ds \partial_{q_i} \tilde{f}_{q_i} =  -\sum_{j\neq i} \partial_{q_j} \tilde{f}_{q_i}.
    \end{cases}
\end{equation*}

In order to keep the notations as light as possible, we will write $\Diff f$ in place of $\Diff \tilde{f}$.

\subsection{Defining the exterior fields}

We recall that in order to apply Theorem~\ref{prop:argFinal}, we need to show that every blob $\omega_i$ is a solution of a problem \eqref{eq:systemDuPB} with some exterior field $F_i$.

Let $\omega_0$ satisfying the general hypotheses \eqref{hyp:omega0} for $N \ge 3$ and $Z^*$ the stationary vortex crystal presented in Section~\ref{sec:vortex_crystals}. Let $\omega$ be a solution of \eqref{eq:gSQG}. We observe that each blob $\omega_i$ is solution of \eqref{eq:gSQG} by letting
\begin{equation*}
    F(x,t) = F_i(x,t) = \int \sum_{j \neq i} K_\alpha(x,y)\omega_j(y,t)\dd y.
\end{equation*}
Since $\nabla_x\cdot K_\alpha(x,y) = 0$, then $\nabla\cdot F =0$. Moreover, we have the following lemma, that proves that Hypothesis $(ii)$ of Theorem~\ref{prop:argFinal} is satisfied.
\begin{lemme}
Let $i \in \{1,\ldots,N\}$. We have for all $t \le \tau_{\eps,\beta}$ that
\begin{equation*}
    \big| F_i(B_i(t),t) - f_i(B(t)) \big| \le C \sum_{j=1}^N\sqrt{I_j},
\end{equation*}
where $C$ depends only on $\alpha$, the $a_i$ and $Z^*$.
\end{lemme}
\begin{proof}
\begin{align*}
    \left|F_i(B_i(t),t)- f_i(B(t))\right|  & = \left|\sum_{j\neq i} \int K_\alpha(B_i(t),y) \omega_j(y,t)\dd y -  \sum_{j\neq i} a_j K_\alpha(B_i(t),B_j(t))\right|\\
    & = \left|\sum_{j\neq i}\int \big( K_\alpha(B_i(t),y) - K_\alpha(B_i(t),B_j(t)) \big)\omega_j(y,t)\dd y\right| \\
    & \le \sum_{j\neq i} \int \big| K_\alpha(B_i(t),y) - K_\alpha(B_i(t),B_j(t)) \big|\frac{|a_j|\omega_j(y,t)}{a_j}\dd y \\
    & \le C \sum_{j\neq i}|a_j|\int \big| y-B_j(t) \big|\frac{\omega_j(y,t)}{a_j}\dd y \\
    & \le C \sum_{j\neq i} \sqrt{I_j(t)},
\end{align*}
where we used that $K_\alpha$ and its derivatives are smooth on $\R^2\times\R^2 \setminus \{ x = y \}$, and go to 0 when $|x-y| \to +\infty$, so in particular $K_\alpha$ is a Lipschitz map on $\bigcup_{i \neq j} D(z_i^*,\eps^\beta)\times D(z_j^*,\eps^\beta)$.
\end{proof}
The existence of $\kappa_0$, $\kappa_1$ and $\kappa_2$ is trivial since $F$ and $f$ are Lipschitz maps while the blobs and point-vortices remain far from each other, which is always the case when $t \le \tau_{\eps,\beta}$ and $t\le \tau_Z$, so Hypothesis $(iii)$ is always satisfied.

Therefore, in order to apply Theorem~\ref{prop:argFinal}, we come down to prove that there is indeed an eigenvalue of $\Diff f(Z^*)$ with positive real part. Then the only remaining difficult task is to estimate the constants $\kappa_1$, $\kappa_2$ and $\lambda_0$ to obtain an information on the lowest possible choice of $\nu$.

\subsection{Proof of Theorems~\ref{theo:N=3} and~\ref{theo:N=3sqg}}

We now prove Proposition~\ref{prop:N=3}, which in turns proves Theorem~\ref{theo:N=3}. We construct explicitly the vortex crystal as described in Section~\ref{sec:vortex_crystals} with $N=3$, namely the configuration:
\begin{equation*}
    \begin{cases}
        z_1^* = (-1,0), &  a_1 = 1 \\
        z_2^* = (1,0), & a_2 = 1\\
        z_3^* = (0,0), & a_3 = -\frac{1}{2^\alpha}.
    \end{cases}
\end{equation*}

We now compute $\Diff f(Z^*)$ using the method previously described at the end of Section~\ref{sec:vortex_crystals} to obtain the $6\times 6$ matrix:

\begin{equation*}
   \Diff f(Z^*) =  C_\alpha
\left(
\begin{array}{cccccc}
 0 & 2^{-(\alpha+1)} & 0 & 2^{-(\alpha+1)} & 0 & -2^{-\alpha } \\
 2^{-(\alpha+1)} \alpha  & 0 & 2^{-(\alpha+1)} \alpha  & 0 & -2^{-\alpha } \alpha  & 0 \\
 0 & 2^{-(\alpha+1)} & 0 & 2^{-(\alpha+1)} & 0 & -2^{-\alpha } \\
 2^{-(\alpha+1)} \alpha  & 0 & 2^{-(\alpha+1)} \alpha  & 0 & -2^{-\alpha } \alpha  & 0 \\
 0 & 1 & 0 & 1 & 0 & -2 \\
 \alpha  & 0 & \alpha  & 0 & -2 \alpha  & 0 \\
\end{array}
\right)
\end{equation*}
The eigenvalues of this matrix are 0, with multiplicity 4, and $\pm  C_\alpha(2-2^{-\alpha} )\sqrt{\alpha}$. So by letting $\lambda_0 = C_\alpha(2-2^{-\alpha} )\sqrt{\alpha}$ we have here a positive eigenvalue, associated with the eigenvector
\begin{equation*}
    v_{\lambda_0} = \left(-1,\sqrt{\alpha},-1,\sqrt{\alpha},-2^{\alpha+1},\sqrt{\alpha}2^{\alpha+1}\right).
\end{equation*}
This conclude Step 1.

We now apply Theorem~\eqref{prop:argFinal} and obtain the existence for any $\alpha \in [1,2)$ of $\nu \ge 2$ and $\beta_0< \frac{4-2\alpha}{5-\alpha}$ such that for every $\beta \in (\beta_0,1)$, every
\begin{equation*}
    \xi_1 > \frac{1-\beta}{\lambda_0} = \frac{1-\beta}{C_\alpha(2-2^{-\alpha} )\sqrt{\alpha}}
\end{equation*}
and every $\eps >0$ small enough, there exists $\omega_0$ satisfying~\eqref{hyp:omega0} such that any solution $\omega$ of \eqref{eq:gSQG} satisfies
\begin{equation*}
    \tau_{\eps,\beta} \le \xi_1 |\ln\eps|.
\end{equation*}
We proved Theorem~\ref{theo:N=3sqg}, and by taking $\alpha=1$ we proved Theorem~\ref{theo:N=3} as well.


\section{Unstable point-vortex in a bounded domain}\label{sec:proofbounded}
In this section we construct a bounded domain $\Omega$ and an initial data $\omega_0$ satisfying \eqref{hyp:omega0} with $N=1$, such that the solution of \eqref{eq:Eu} satisfies that
\begin{equation*}
    \tau_{\eps,\beta} \le \xi_1 |\ln \eps|,
\end{equation*}
for $\eps$ small enough and for some $\xi_1$. Since $N=1$, we denote by $z^* = z_1^* = Z^*$ and $a=a_1$.

We start by recalling some facts about the point-vortex dynamics in bounded domain, then we construct the domain. Finally, we prove Theorem~\ref{theo:boundedDomain} using again Theorem~\ref{prop:argFinal}.

\subsection{Euler equations and point-vortices in bounded domains}

For the rest of the paper, we consider the Euler equations~\eqref{eq:Eu}, which differs from before in the sense that now $\alpha= s = 1$, and $ \omega : \Omega \to \R$, where $\Omega$ is a bounded simply connected subset of $\R^2$. We now recall that the problem
\begin{equation*}
\begin{cases}
    \Delta\Psi = \omega & \quad \text{on } \Omega \\
    \Psi = 0 & \quad \text{on } \partial\Omega
\end{cases}
\end{equation*}
has a unique solution
\begin{equation*}
    \Psi(x) = \int_\Omega G_\Omega(x,y)\omega(y)\dd y,
\end{equation*}
where $G_\Omega : \Omega\times\Omega \setminus \{x=y\} \to \R$ is the Green's function of $-\Delta$ with Dirichlet condition in the domain $\Omega$.
Therefore, we have the Biot-Savart law:
\begin{equation*}
    u(x,t) = \int_\Omega \nabla_x^\perp G_\Omega(x,y)\omega(y,t)\dd y.
\end{equation*}
An important property of the Green's function $G_\Omega$ is that it decomposes as
\begin{equation*}
    G_\Omega = G_{\R^2} + \gamma_\Omega,
\end{equation*}
where 
\begin{equation*}
    G_{\R^2}(x,y) = \frac{1}{2\pi}\ln|x-y|,
\end{equation*}
and $\gamma_\Omega : \Omega\times\Omega \to \R_+$ harmonic in both variables. We denote by $\tgm_\Omega : x \mapsto \gamma_\Omega(x,x)$ the Robin's function of the domain $\Omega$. This map plays a crucial role in the study of the point-vortex dynamics in bounded domain. Indeed, the point-vortices in move according to the system
\begin{equation*}
    \forall i \in \llbracket 1,N\rrbracket, \quad \der{}{t} z_i(t) = \sum_{\substack{1 \le j \le N \\ j\neq i}}  a_j \frac{\big(z_i(t)-z_j(t)\big)^\perp}{\big| z_i(t) - z_j(t)\big|^2} + \sum_{j=1}^N a_j \nabla^\perp_x\gamma_\Omega\big(z_i(t),z_j(t)\big),
\end{equation*}
which in the case $N=1$ reduces to
\begin{equation*}
    \der{}{t} z(t) = \frac{a}{2} \nabla^\perp \tgm_\Omega(z(t)).
\end{equation*}

Our plan to prove Theorem~\ref{theo:boundedDomain} is the same as the proof of Theorem~\ref{theo:N=3}. The aim is to apply Proposition\ref{prop:solution_divergente} and Theorem~\ref{prop:argFinal}. We then straight away notice that when $\omega_0$ satisfy \eqref{hyp:omega0} with $N=1$, and thus $\omega$ (which we can extend by 0 to $\R^2$) is constituted of a unique blob that solves \eqref{eq:systemDuPB} by setting
\begin{equation*}
    F(x,t) = \int \nabla_x^\perp\gamma_\Omega(x,y) \omega(y,t) \dd y,
\end{equation*}
and $z$ is a solution of \eqref{eq:pbz} by setting
\begin{equation*}
    f(z) = \frac{a}{2} \nabla^\perp\tgm_\Omega(z(t)).
\end{equation*}
Therefore, we are looking for $z^*$ an unstable critical point of the Robin's function $\tgm_\Omega$, namely such that $\nabla\tgm_\Omega(z^*) = 0$ and $\Diff\nabla^\perp \tgm_\Omega(z^*)$ has a positive eigenvalue. However the Robin's function is not known explicitly in general, and the existence of such a point depends on the domain $\Omega$. Fortunately, we recall that for any simply connected domain $\Omega$ and for any $z^* \in \Omega$, there exists a biholomorphic map $T : \Omega \to D:= D(0,1)$ such that $T(z_0) = 0$. Such a map also satisfy that
\begin{equation*}
    \forall x \neq y \in \Omega, \quad G_\Omega(x,y) = G_D\big(T(x),T(y)\big),
\end{equation*}
and thus 
\begin{equation*}
    \tgm_\Omega(x) = \tgm_D(x) + \frac{1}{2\pi} \ln|T'(x)|.
\end{equation*}

\subsection{Known results on critical points of the Robin's function}\label{sec:donati_iftimie}
In this section we refer to \cite{Donati_Iftimie_2020} and recall the following results. For more details on the Robin's function we refer the reader to \cite{Gustafsson_1979,flucher_variational_1999}.
\begin{prop}[\cite{Donati_Iftimie_2020}, Proposition 2.4]
    If $T:\Omega \to D$ is a biholomorphic map such that $T(z^*) = 0$, then $$\nabla\tgm_\Omega(z^*) = 0 \Longleftrightarrow T''(z^*) = 0.$$
\end{prop}
In particular, if the domain $\Omega$ has two axes of symmetry, then the intersection $z^*$ necessarily satisfy $\nabla\tgm_\Omega(z^*) = 0$ and thus $T''(z^*) = 0$. 
For the time being, we assume the existence of such $z^*$ and state some of its properties.
\begin{lemme}[\cite{Donati_Iftimie_2020}, Lemma 4.1]
Let $T : \Omega \to D(0,1)$ and $z^* \in \Omega$ such that $T(z^*)=0$ and $T''(z^*) = 0$, then for every $\eps > 0$ small enough,
\begin{equation*}
    \forall x,y,z \in D(z^*,\eps^\beta), \quad \left|\nabla_x^\perp\gamma_\Omega(x,y) -\nabla_x^\perp\gamma_\Omega(z,y) \right| = |x-z|\left(\frac{|T'''(z^*)|}{6\pi |T'(z^*)|} + \mathcal{O}\left(\eps^\beta\right)\right).
\end{equation*}
\end{lemme}
The direct corollary of this lemma is that $F$ satisfies \eqref{eq:Lip_F(kappa1)} with $$\kappa_1 = |a|\frac{|T'''(z^*)|}{6\pi |T'(z^*)|} + o(1).$$
\begin{prop}[\cite{Donati_Iftimie_2020}, Sections 3.2 and 3.3]
    Let $\Omega \subset \R^2$ be a bounded simply connected domain. Then for any biholomorphism $T : \Omega \to D$ mapping $z^*$ to $0$ such that $T''(z^*) = 0$, the hessian matrix $\Diff^2\tgm_\Omega(z^*)$ has non degenerate eigenvalues of opposite signs if and only if $|T'''(z^*)| > 2|T'(z^*)|^3$. In that case, these eigenvalues are
\begin{equation*}
    \lambda_\pm = \frac{2|T'(z^*)|^2 \pm |T'''(z^*)|/|T'(z^*)|}{2\pi},
\end{equation*}
 and the eigenvalues of $\Diff\nabla^\perp\tgm_\Omega(z^*)$ are $\pm \sqrt{-\lambda_+\lambda_-}$, so that 
 \begin{equation*}
     \lambda_0 =  \frac{|a|}{4\pi |T'(z^*)|}\sqrt{|T'''(z^*)|^2 -4|T'(z^*)|^6}
 \end{equation*}  
 is a positive eigenvalue of $\Diff f(z^*)$.
\end{prop}
From this we deduce two things. First, the condition $|T'''(z^*)| > 2|T'(z^*)|^3$ is a criteria to establish that $\Diff f(z^*)$ has an eigenvalue with positive real part. Second, since $\Diff^2 \tgm_\Omega$ is a real symmetric matrix, and since $\Diff \nabla^\perp(\tgm_\Omega) = R_{\pi/2}\Diff^2 \tgm_\Omega$ with $R_{\pi/2} = \begin{pmatrix} 0 & -1 \\ 1 & 0 \end{pmatrix}$, then
\begin{equation*}
    \ntriple{\Diff f(z^*)} = \frac{|a|}{2}\ntriple{\Diff^2 \tgm_\Omega} = \frac{|a|}{2}\lambda_+,
\end{equation*}
so that $f$ satisfies \eqref{eq:Lip_f(kappa2)} with
\begin{equation*}
    \kappa_2 = |a|\frac{2|T'(z^*)|^2 + |T'''(z^*)|/|T'(z^*)|}{4\pi}.
\end{equation*}
Therefore, in view of relation~\eqref{eq:nu}, we will be able in our construction to choose any $\nu$ such that
\begin{equation*}
    \nu > \frac{\frac{5}{3}|T'''(z^*)| + 2|T'(z^*)|^3}{\sqrt{|T'''(z^*)|^2 -4|T'(z^*)|^6}} + 1.
\end{equation*}
which satisfies in particular that
\begin{equation*}
    \frac{|T'''(z^*)|}{|T'(z^*)|^3} > \frac{15 + 9\sqrt{65}}{28} \; \Longrightarrow \; \frac{\frac{5}{3}|T'''(z^*)| + 2|T'(z^*)|^3}{\sqrt{|T'''(z^*)|^2 -4|T'(z^*)|^6}} + 1 < 4.
\end{equation*}
In conclusion of this section, in order to prove Theorem~\ref{theo:boundedDomain}, we need in particular to construct a domain $\Omega$ satisfying the existence of a point $z^*$ and a biholomorphic map $T : \Omega \to D$ such that
\begin{equation*}
    T(z^*) = T''(z^*) = 0 \; \text{ and } \; \frac{|T'''(z^*)|}{|T'(z^*)|^3} > 2
\end{equation*}
for the construction to be possible with some $ \nu \ge 2$, and that
\begin{equation*}
     \frac{|T'''(z^*)|}{|T'(z^*)|^3} > \frac{15 + 9\sqrt{65}}{28} := c_0 \approx 3.12,
\end{equation*}
for the construction to be possible with $\nu =4$.

\subsection{Construction of the domain}\label{sec:domain}
\begin{figure}
    \centering
    \includegraphics[width=0.5\textwidth]{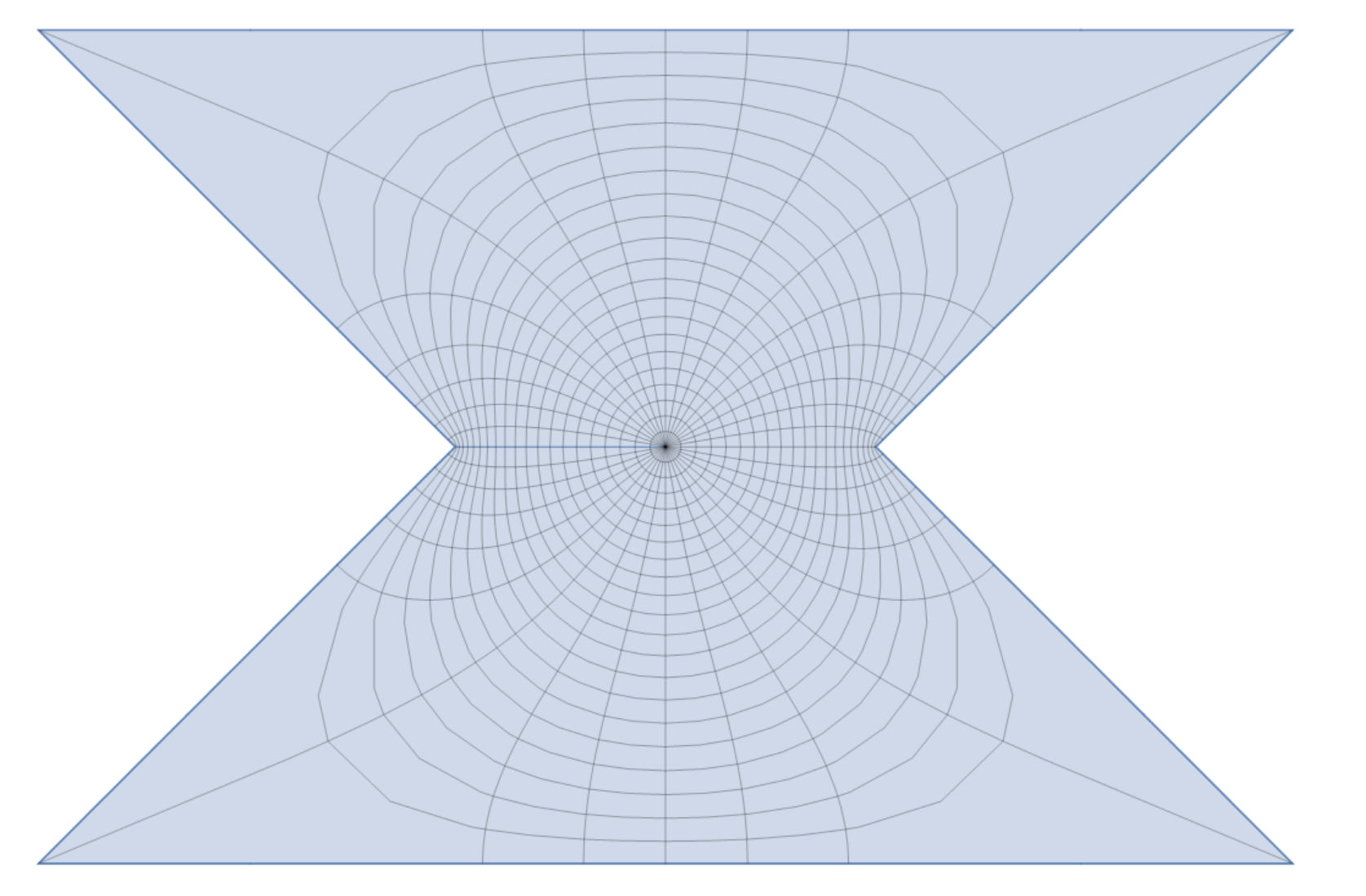}
    \caption{The domain $\Omega_{\frac{3}{4}}$.}
    \label{fig:domain}
\end{figure}

Let $\delta \in [\frac{1}{2},1)$. Let $\zeta = e^{i\frac{\pi}{3}}$. Using the Schwartz Christoffel formula (see for instance \cite{Driscoll_Trefethen_2002_Schwarz-Christoffel_Mapping}), we define the conformal map $S_\delta : D \to \Omega_\delta := S_\delta(D)$ mapping $0$ to $0$ such that
\begin{align*}
    S_\delta'(z) & := \frac{1}{(z-1)^{(1-2\delta)}(z-\zeta)^\delta(z-\zeta^2)^\delta(z+1)^{(1-2\delta)}(z-\zeta^4)^\delta(z-\zeta^5)^\delta} \\
    & = \frac{(z^2-1)^{3\delta -1}}{(z^6-1)^\delta}.
\end{align*}
Let $T_\delta = S_\delta^{-1}$. We compute that $|T_\delta'(0) = |S_\delta'(0)| = 1$, $T_\delta''(0) = S_\delta''(0) = 0$ and $$|T_\delta'''(0)| = |S_\delta'''(0)| = 6\delta-2.$$ Therefore, we have first of all that \begin{equation*}
    |T_\delta'''(0)|  > 2 |T_\delta'(0)|^3 \Longleftrightarrow \delta > \frac{2}{3}.
\end{equation*}
So the domain $\Omega_\delta$ satisfies that $\Diff f(0)$ has an eigenvalue with positive real part and its Robin's function has a saddle point. In Figure~\ref{fig:domain} is plotted the domain $\Omega_{\frac{3}{4}}$. We can then do our construction of $\omega_0$ in $\Omega_{\frac{3}{4}}$, but it's not enough to it with $\nu =4$.

However, we have that
\begin{equation*}
    \frac{|T_\delta'''(0)|}{|T_\delta'(0)|^3} = 6\delta -2.
\end{equation*}
In particular we observe that for $\delta \in \Big( \frac{c_0+2}{6} , 1\Big)$, $\nu = 4$ satisfies \eqref{eq:nu}.
Since $0.9 > \frac{c_0+2}{6} $, the domain $\Omega_{0.9}$ is a suitable domain to do the construction of $\omega_0$ with $\nu =4$. More details and illustrations about the family of domains $(\Omega_\delta)_\delta$ are given in Appendix~\ref{sec:appendix_domain}.

To obtain a smooth domain $\tilde{\Omega}_d$ with the same properties, one takes an increasing sequence $\Omega_n$ of smooth domains which are symmetric with respect to both axes and converge towards $\Omega_{9/10}$. By symmetry, 0 is necessarily a critical point of the Robin's function of every domain. Then, we introduce $T_n$ the sequence of conformal maps mapping $\Omega_n$ to $D(0,1)$ satisfying $T_n(0) = 0$ and $T_n'(0) \in \R_+$. The construction can be done so that $T_n \to T_{9/10}$ locally in every $C^k$, $k \ge 0$, so that
\begin{equation*}
    \frac{|T_n'''(0)|}{|T_n'(0)|^3} \tend{n}{+\infty} \frac{|T_{9/10}'''(0)|}{|T_n'(0)|^3} > c_0,
\end{equation*}
so there exists $n_0 \in \N$, such that the smooth domain $\Omega_{n_0}$ is such that
\begin{equation*}
        4 > \frac{\frac{5}{3}|T_{n_0}'''(0)| + 2|T_{n_0}'(0)|^3}{\sqrt{|T_{n_0}'''(0)|^2 -4|T_{n_0}'(0)|^6}} + 1.
\end{equation*}

\subsection{Proof of Theorem~\ref{theo:boundedDomain}}
Let $\delta \in \Big( \frac{c_0+2}{6} , 1\Big)$, let $\Omega := \Omega_\delta$ (or $\Omega_{n_0}$ as described in the previous section to work with a smooth domain), $z^* = 0$ and $a=1$.

From Sections~\ref{sec:donati_iftimie} and~\ref{sec:domain}, $\Diff f(0)$ has an eigenvalue $\lambda_0 = \frac{1}{4\pi} \sqrt{ 12\delta(3\delta-2)} > 0,$ so Hypothesis $(i)$ of Theorem~\ref{prop:argFinal} is satisfied.
We recall that Hypothesis $(iii)$ is satisfied since $f$ and $F$ are Lipschitz maps far from the boundary.

Therefore, there only remains to prove that Hypothesis $(ii)$ is satisfied to apply Theorem~\ref{prop:argFinal} and conclude the proof of Theorem~\ref{theo:boundedDomain}.

Let us compute.
\begin{align*}
    \big| F(B(t),t) - f(B(t)) \big| & = \left| \int \nabla_x^\perp\gamma_\Omega(B(t),y)\omega(y,t) \dd y - \frac{a}{2} \tgm_\Omega(B(t))\right| \\
    & = \left| \int \big(\nabla_x^\perp\gamma_\Omega(B(t),y) - (\nabla_x^\perp\gamma_\Omega(B(t),B(t))\big)\omega(y,t)\dd y\right| \\
    & \le C|a|\int |y - B(t)|\frac{\omega(y,t)}{a}\dd y\\
    & \le C \sqrt{I(t)},
\end{align*}
where on the last line the constant $C$ depends only on $\Omega$. 
Letting $\delta \to 1$, we observe that $\lambda_0 \to \frac{\sqrt{3}}{2\pi}$. Now applying Theorem~\eqref{prop:argFinal}, we have that for any $\beta_0 < 1/2$, one can construct a suitable $\omega_0$ satisfying \eqref{hyp:omega0} with $N=1$ and $\nu =4$. Theorem~\ref{theo:boundedDomain} is now completely proved.

\appendix

\section{The family of biconvex hexagonal domains and their Robin's function}\label{sec:appendix_domain}
Let us mention that such domains were drawn and studied already in \cite{Flucher_1992_Extremal_function_for_Trudinger_Moser}.

\begin{figure}[H]
    \centering
    \includegraphics[width=0.3\textwidth]{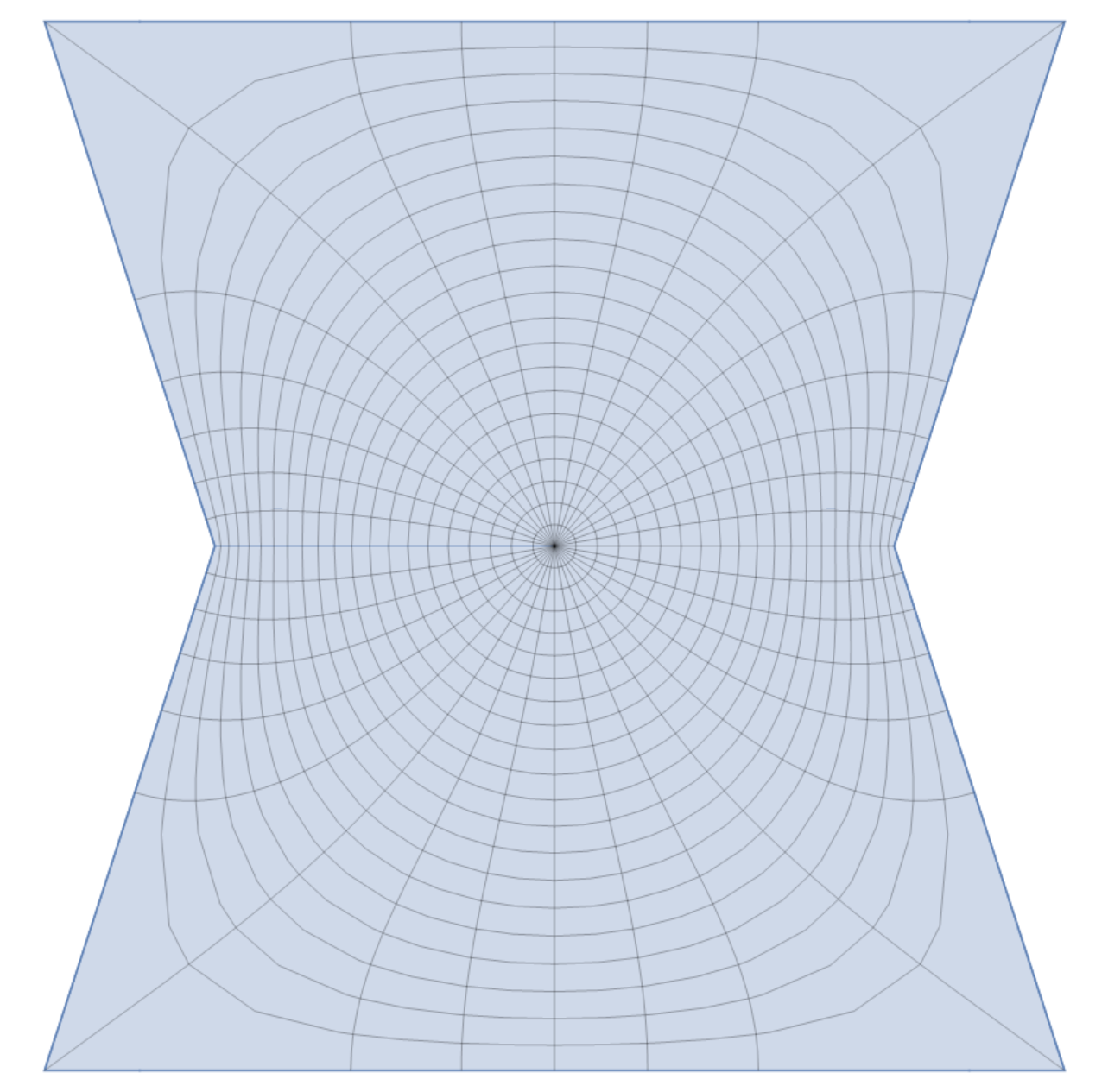}
    \includegraphics[width=0.3\textwidth]{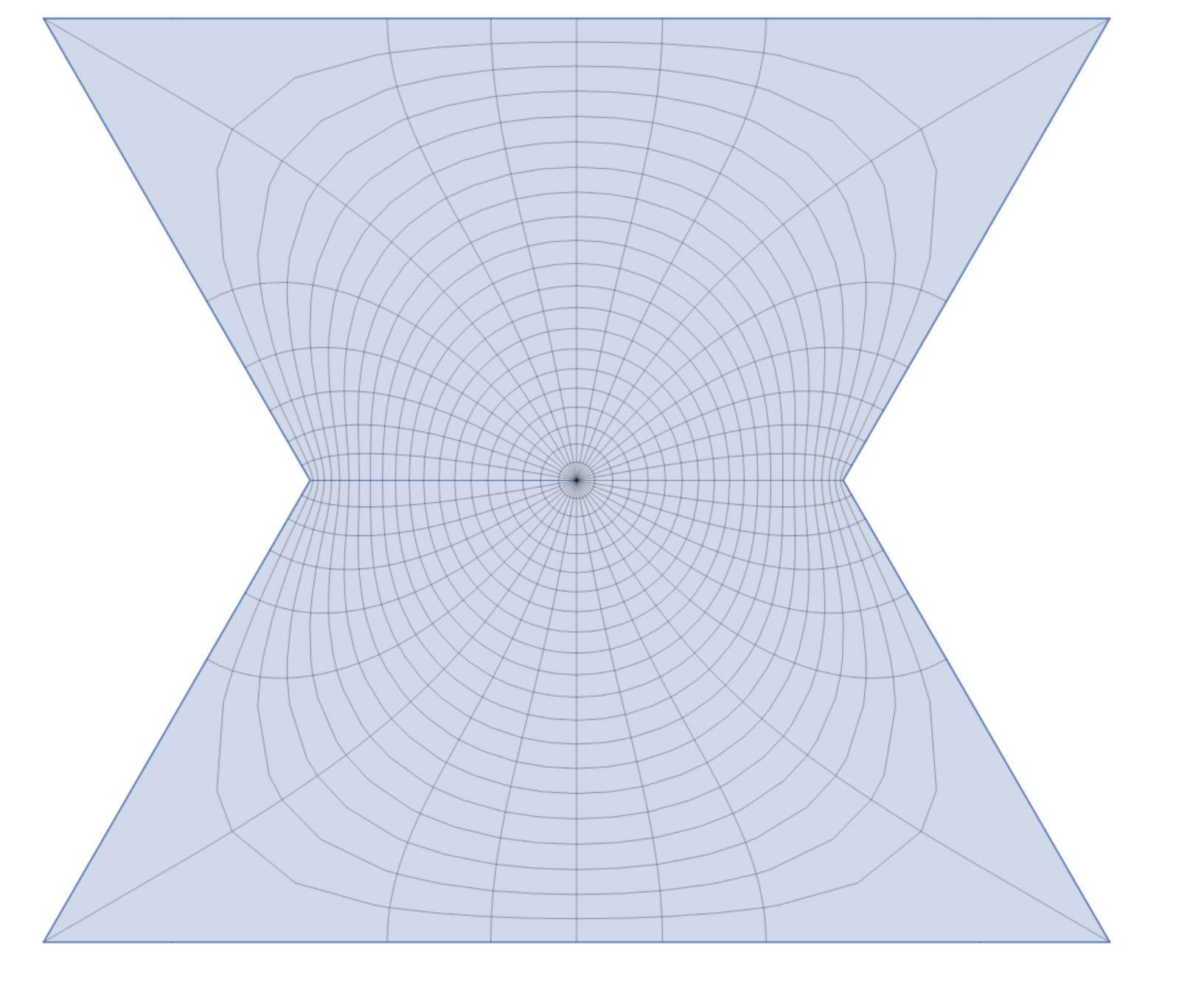}
    \includegraphics[width=0.3\textwidth]{Figures/Domain_3-4.PNG}
    \caption{Plot of the domain $S_\delta(D)$, for, left to right: $\delta = \frac{2}{5}$, $\delta = \frac{2}{3}$ and $\delta = \frac{3}{4}$.}
    \label{fig:domains}
\end{figure}

\begin{figure}[H]
    \centering
    \includegraphics[width=0.3\textwidth]{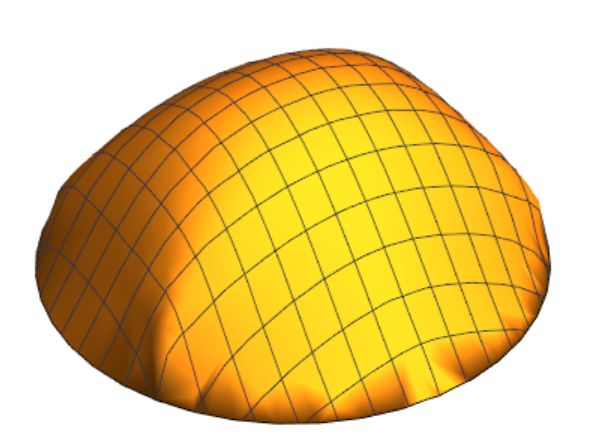}
    \includegraphics[width=0.3\textwidth]{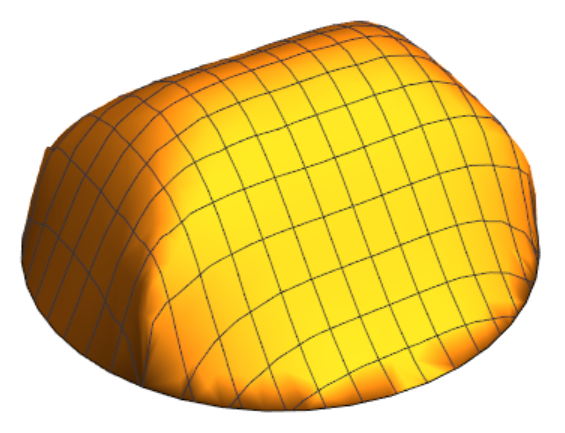}
        \includegraphics[width=0.3\textwidth]{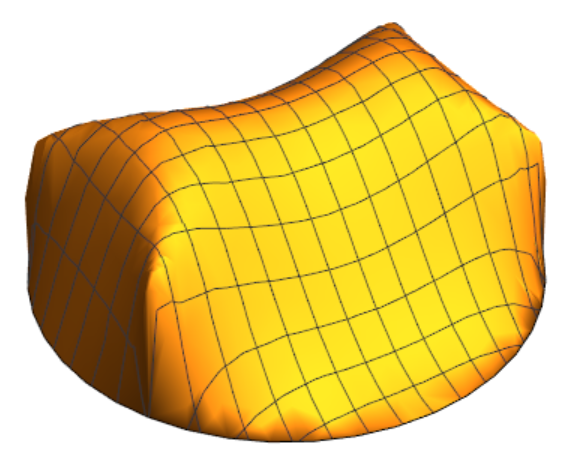}
    \caption{Plot of the map $r_{S_\delta(D)} \circ S_\delta : x \mapsto |S_\delta'(x)|(1-|x|^2)$, for, left to right: a stable case ($\delta = \frac{2}{5}$), the critical case ($\delta = \frac{2}{3})$, and an unstable case ($\delta = \frac{9}{10}$). }
    \label{fig:robin}
\end{figure}

We used Wolfram Mathematica to plot several domains in Figure \ref{fig:domains}. In those cases, the Robin's function is not a very nice function to plot. Instead we introduce the conformal radius:
\begin{equation*}
    r_\Omega(x) = e^{-2\pi\tgm_\Omega(x)}.
\end{equation*}
It satisfies the transfer formula (see \cite{Flucher_1992_Extremal_function_for_Trudinger_Moser}):
\begin{equation*}
    r_\Omega(S(x)) = |S'(x)|r_D(x) = |S'(x)|(1-|x|^2).
\end{equation*}
This map is a lot easier to draw, see Figure \ref{fig:robin}. We see that the Robin's function of the domains obtained with $\delta > 2/3$ have a saddle point in 0.

\section{Technical lemmas}\label{sec:appendix_proofs}

\subsection{Actual construction of the initial data}
\text{ }

We formulate the following remark.
\begin{remark}\label{prop:choixOmega0}
Let $N \ge 1$, $Z^* \in (\R^2)^N$ and $a_i \in \R^*$. For any $Z_0$ such that $|Z_0  -Z^*| = \frac{\eps}{2}$, one can always chose $\omega_0 \in C^\infty(\R^2)$ such that
\begin{itemize}
    \item $\omega_0$ satisfies \eqref{hyp:omega0},
    \item $B(0)=Z_0$,
    \item $\forall i \in \{1,\ldots,N\}, \; I_i(0) \le \eps^{\nu}$.
\end{itemize}
\end{remark}
\begin{proof}
For $N=1$, $z_0 \in \R^2$, $a \in \R^*$, we introduce the vortex patch $\omega_0 = 16a\frac{\eps^{-\nu}}{\pi}\Ind_{D(z_0,\frac{\eps^{\nu/2}}{4})}$ for $\eps$ small enough. We verify that 
\begin{equation*}
    |\omega_0| \le C \eps^{-\nu}
\end{equation*}
with $C = \frac{4a}{\pi}$,
\begin{equation*}
    \int \omega_0(x)\dd x = a,
\end{equation*}
\begin{equation*}
    B(0)= \frac{1}{a}\int_{\R^2} x\omega_0(x)\dd x = z_0,
\end{equation*}
and
\begin{equation*}
    I(0) = \frac{1}{a}\int_{\R^2} |x-B(0)|^2\omega_0(x)\dd x = \int_0^{\frac{\eps^{\nu/2}}{4}}\int_0^{2\pi} 16r^2 \frac{\eps^{-\nu}}{\pi} r\dd \theta \dd r = 8\eps^{-\nu}\frac{\eps^{2\nu}}{256} = \frac{1}{32}\eps^{\nu}.
\end{equation*}
For $N > 1$ we sum patches of this exact form. One can also construct a smooth $\omega_0$ satisfying those constraints, by taking a convenient radially mollified version of these vortex patches such that their support lies within $D(z_i,\eps)$. Since $\frac{\eps^{\nu/2}}{4} < \eps/2$ as $\nu \ge 2$, it is always possible. 
\end{proof}

\subsection{Variant of the Gronwall's inequality}

From the Gronwall's inequality, we can write the following.
\begin{lemme}\label{lemme:gronwall}
Let $f$ be a $C^1$ map and $g$ be positive non decreasing such that 
\begin{equation*}
    |f(x) - f(y)| \le \kappa |x-y|.
\end{equation*}
Let $z$ is a solution of
\begin{equation*}
    z'(t) = f(z(t)),
\end{equation*}
and let $y$ such that
\begin{equation*}
\begin{cases}
    \left|y'(t) - f(y(t))\right| \le g(t) \\
    y(0) = z(0),
\end{cases}
\end{equation*}
where $g : \R_+ \to \R_+$ is smooth.
Then
\begin{equation*}
    |y(t) - z(t)| \le e^{\kappa t}\int_0^t g(s)\dd s .
\end{equation*}
\begin{proof}
On has that 
\begin{align*}
    |y(t)-z(t)| & = \left|\int_0^t \big(y'(s) + z'(s) \big)\dd s \right|
    \\ & \le \int_0^t g(s)\dd s + \left| \int_0^t \big(f(y(s)) - f(z(s))\big) \dd s \right| 
    \\ & \le \int_0^t g(s)\dd s  + \kappa \int_0^t |y(s)-z(s)|\dd s,
\end{align*}
so using now the classical Gronwall's inequality, since $t\mapsto \int_0^t g(s)\dd s$ is non decreasing, we have that
\begin{equation*}
    |y(t)-z(t)| \le e^{\kappa t}\int_0^t g(s) \dd s .
\end{equation*}
\end{proof}
\end{lemme}

\section{Computation of the constants}\label{sec:compConsts}

In Theorem~\ref{prop:argFinal}, we proved that the construction is possible as long as
\begin{equation*}
    \nu > \frac{2}{5-\alpha} \left( (1+\alpha)\frac{\kappa_1+\kappa_2}{\lambda_0} + 4-2\alpha\right),
\end{equation*}
and
\begin{equation*}
    \xi_1 > \frac{1-\beta}{\lambda_0}.
\end{equation*}
Since $\kappa_1$, $\kappa_2$ and $\lambda_0$ ultimately depend on the chosen configuration of point-vortices, so do the bounds on $\nu$ and $\xi_1$.

We now give details and improvements on those bounds.

\subsection{Results}

We state a few results that will be proved in the following sections.

By computing the constants $\kappa_1$ and $\kappa_2$ for the construction done in Section~\ref{sec:proof_SQG} with $N=3$, we obtain the following details on the bound on $\nu$.
\begin{prop}\label{prop:N=3}
    One can achieve the construction of Theorem~\ref{theo:N=3sqg} for any $\alpha \in[1,2)$ with any $\nu$ satisfying
    \begin{equation*}
    \nu > \frac{2}{5-\alpha} \left( (1+\alpha)\frac{2 + 2^{-\alpha}\alpha\sqrt{3(1+2^{1+2\alpha})}}{(2-2^{-\alpha} )\sqrt{\alpha}}+ 4-2\alpha\right),
\end{equation*}
which is greater than $4$.
\end{prop}
\begin{prop}\label{prop:N=7}
Let $\nu =4$. There exists $\alpha_0 >1$ such that for any $\alpha \in [1,\alpha_0)$, there exists an initial configuration $\big((z_i^*)_i,(a_i)_i\big)$ of point-vortices with $N=7$ and $\beta_0 < \frac{4-2\alpha}{5-\alpha}$ such that for every $\beta \in (\beta_0,1)$ there exists $\xi_1$ such that for every $\eps >0$ small enough, there exists $\omega_0$ satisfying \eqref{hyp:omega0} such that any solution $\omega$ of \eqref{eq:gSQG} satisfies
\begin{equation*}
    \tau_{\eps,\beta} \le \xi_1 |\ln \eps|.
\end{equation*}
In particular, if $\alpha=1$, one can take $\xi_1 >  \frac{4\pi}{9}(1-\beta)$.
\end{prop}
Unfortunately, we fail to obtain a rigorous estimate of $\alpha_0$. However, in Section~\ref{sec:N>3}, we numerically check that a construction using 9 blobs can be done with $\nu =4$ for any $\alpha \in [1,2)$. The constant $\nu = 4$ is not optimal.

\subsection{How to compute the constants}

We start by giving a general method on how to obtain $\kappa_1$ and $\kappa_2$. Let $N\ge3$.

We recall that for all $t \le \tau_{\eps,\beta}$, $\supp \omega \subset \bigcup_{j=1}^N D(z_j^*,\eps^\beta)$. Therefore  applying \eqref{eq:developpmentK} to $\tilde{z_i} = x$, $\tilde{x} = x'-x$, $\tilde{z_j} = y$, $\tilde{y} = 0$, we have $\forall x,x' \in D(z^*,\eps^\beta)$, and $\forall t \le \tau_{\eps,\beta}$,
\begin{align*}
    F_i(x,t)-F_i(x',t) & = \int \sum_{j\neq i}\big(K_\alpha(x,y) - K_\alpha(x',y)\big) \omega_j(y,t) \dd y
    \\  & = \sum_{j\neq i}C_\alpha\int\frac{(x-y)^\perp}{|x-y|^{\alpha+3}}(\alpha+1)(x'-x)\cdot(x-y)\omega(y,t) \dd y + o(|x-x'|)
\end{align*}
so that
\begin{align*}
     \Big|(x-x')\cdot\big(F_i(x,t)-F_i(x',t)\big)\Big| & \le|x-x'|^2\left(\sum_{j\neq i}\int \frac{C_\alpha}{|x-y|^{\alpha+1}}|\omega(y,t)| \dd y + o(1)\right) \\
     & \le |x-x'|^2\left(\sum_{j\neq i}\frac{C_\alpha|a_j|}{|z_i^*-z_j^*|^{\alpha+1}} + o(1)\right).
\end{align*}
Therefore we obtain that for every $i \in \{1,\ldots,N\}$, for every $x,x' \in D(z^*,\eps^\beta)$ and for every $t \le \tau_{\eps,\beta}$,
\begin{equation*}
    \Big|(x-x')\cdot\big(F_i(x,t)-F_i(x',t)\big)\Big| \le \kappa_1 |x-x'|^2
\end{equation*}
with
\begin{equation*}
    \kappa_1 = C_\alpha\max_{1\le i \le N}\sum_{j\neq i}\frac{|a_j|}{|z_i^*-z_j^*|^{\alpha+1}} + o(1).
\end{equation*}

For $\kappa_2$, we need a Lipschitz type estimate for $f$ on $\bigcup_{j=1}^N D(z_j^*,2\eps^\beta)$, and thus have that
\begin{equation*}
    \kappa_2 = \Big\| \ntriple{\Diff f(Z)} \Big\|_{\infty,\bigcup_{j=1}^N D(z_j^*,2\eps^\beta)} =  \ntriple{\Diff f(Z^*)} + o(\eps^\beta) = \sqrt{\rho\big( \Diff f(Z^*) [\Diff f(Z^*)]^t\big)} + o(\eps^\beta)
\end{equation*}
where $\rho$ is the spectral radius, that is in our case the greatest eigenvalue in absolute value of the real symmetric matrix $\Diff f(Z^*) [\Diff f(Z^*)]^t$.

\subsection{Proof of Proposition~\ref{prop:N=3}}
We now take again $N=3$ and compute.
First, we directly have that
\begin{equation*}
    \kappa_1 = 2C_\alpha +o(1).
\end{equation*}

\begin{figure}[H]
    \centering
    \includegraphics[width=0.5\textwidth]{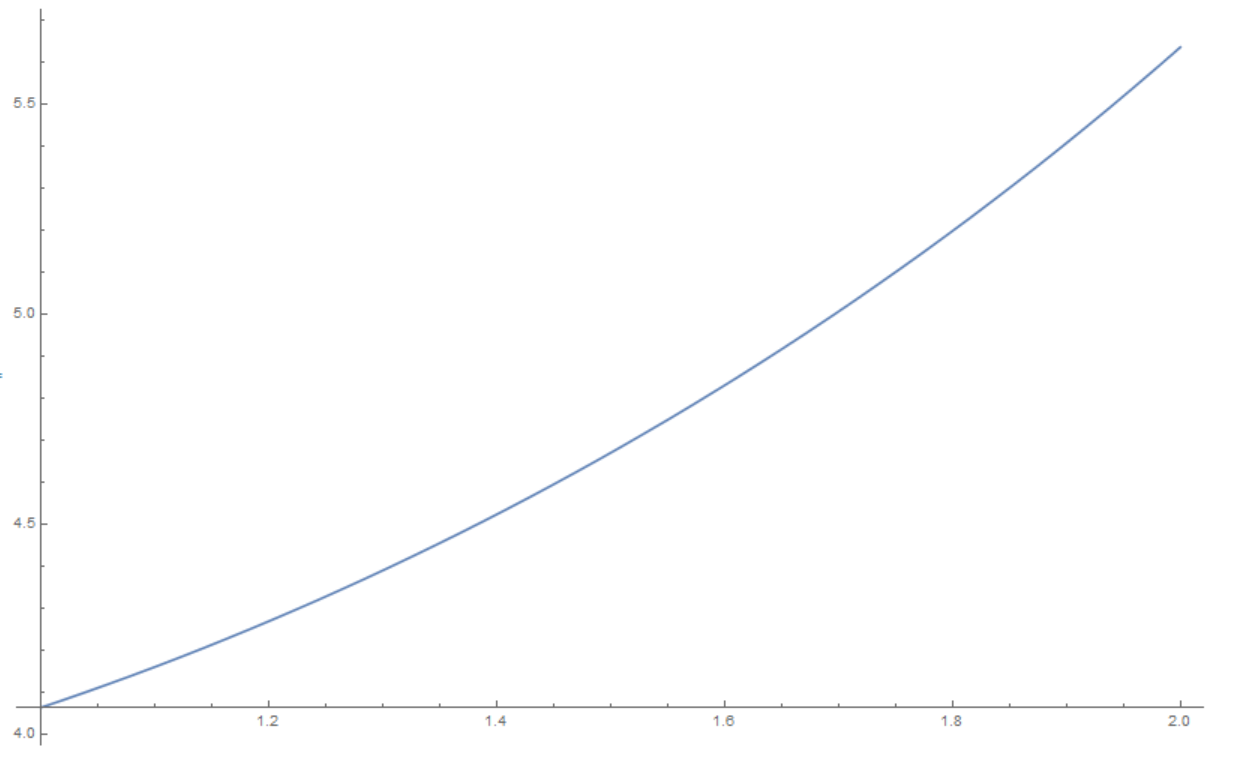}
    \caption{Plot of $g(\alpha)$ in the range $\alpha \in [1,2)$.}
    \label{fig:plot_nu_N=3}
\end{figure}

We now compute the eigenvalues of $\Diff f(Z^*)[\Diff f(Z^*)]^t$ and observe that
\begin{equation*}
    \ntriple{\Diff f(Z^*)} = C_\alpha 2^{-\alpha}\alpha\sqrt{3(1+2^{1+2\alpha})}.
\end{equation*}
Therefore, it is possible to choose $\nu$ such that
\begin{equation*}
    \nu > \frac{2}{5-\alpha} \left( (1+\alpha)\frac{\kappa_1+\kappa_2}{\lambda_0} + 4-2\alpha\right)
\end{equation*}
for $\eps$ small enough as soon as
\begin{equation*}
    \nu > \frac{2}{5-\alpha} \left( (1+\alpha)\frac{2 + 2^{-\alpha}\alpha\sqrt{3(1+2^{1+2\alpha})}}{(2-2^{-\alpha} )\sqrt{\alpha}}+ 4-2\alpha\right) := g(\alpha).
\end{equation*}
The plot of $\alpha \mapsto g(\alpha)$ is given in Figure~\ref{fig:plot_nu_N=3}. 
This concludes the proof of Proposition~\ref{prop:N=3}.

\subsection{Proof of Proposition~\ref{prop:N=7}}\label{sec:N>3}

Using the exact same method, we now construct the vortex crystal with $N=7$.

Since $\forall i\neq j$, $|z_i^* - z_j^*|\ge 1$, then we obtain directly that one can take $\kappa_1 = \frac{3}{\pi} + o(1)$. However in general, we are not able to compute $\kappa_2$ and $\nu$.

We now assume that $\alpha =1$. We then compute that
\begin{multline*}
\Diff f(Z^*) =  \frac{1}{2\pi} \times \\ 
    {\footnotesize\left(
\begin{array}{cccccccccccccc}
 c_1 & -\frac{35}{24} & 0 & 1 & -\frac{1}{2 \sqrt{3}} & \frac{1}{6} & -\frac{\sqrt{3}}{8} & -\frac{1}{8} & 0 & -\frac{1}{3} & \frac{\sqrt{3}}{2} & -\frac{1}{2} & \frac{5 \sqrt{3}}{4} & \frac{5}{4} \\
 -\frac{35}{24} &-c_1 & 1 & 0 & \frac{1}{6} & \frac{1}{2 \sqrt{3}} & -\frac{1}{8} & \frac{\sqrt{3}}{8} & -\frac{1}{3} & 0 & -\frac{1}{2} & -\frac{\sqrt{3}}{2} & \frac{5}{4} & -\frac{5 \sqrt{3}}{4} \\
 0 & 1 &-c_1 & -\frac{35}{24} & -\frac{\sqrt{3}}{2} & -\frac{1}{2} & 0 & -\frac{1}{3} & \frac{\sqrt{3}}{8} & -\frac{1}{8} & \frac{1}{2 \sqrt{3}} & \frac{1}{6} & -\frac{5 \sqrt{3}}{4} & \frac{5}{4} \\
 1 & 0 & -\frac{35}{24} & c_1 & -\frac{1}{2} & \frac{\sqrt{3}}{2} & -\frac{1}{3} & 0 & -\frac{1}{8} & -\frac{\sqrt{3}}{8} & \frac{1}{6} & -\frac{1}{2 \sqrt{3}} & \frac{5}{4} & \frac{5 \sqrt{3}}{4} \\
 -\frac{1}{2 \sqrt{3}} & \frac{1}{6} & -\frac{\sqrt{3}}{2} & -\frac{1}{2} & 0 & \frac{35}{12} & \frac{\sqrt{3}}{2} & -\frac{1}{2} & \frac{1}{2 \sqrt{3}} & \frac{1}{6} & 0 & \frac{1}{4} & 0 & -\frac{5}{2} \\
 \frac{1}{6} & \frac{1}{2 \sqrt{3}} & -\frac{1}{2} & \frac{\sqrt{3}}{2} & \frac{35}{12} & 0 & -\frac{1}{2} & -\frac{\sqrt{3}}{2} & \frac{1}{6} & -\frac{1}{2 \sqrt{3}} & \frac{1}{4} & 0 & -\frac{5}{2} & 0 \\
 -\frac{\sqrt{3}}{8} & -\frac{1}{8} & 0 & -\frac{1}{3} & \frac{\sqrt{3}}{2} & -\frac{1}{2} & c_1 & -\frac{35}{24} & 0 & 1 & -\frac{1}{2 \sqrt{3}} & \frac{1}{6} & \frac{5 \sqrt{3}}{4} & \frac{5}{4} \\
 -\frac{1}{8} & \frac{\sqrt{3}}{8} & -\frac{1}{3} & 0 & -\frac{1}{2} & -\frac{\sqrt{3}}{2} & -\frac{35}{24} &-c_1 & 1 & 0 & \frac{1}{6} & \frac{1}{2 \sqrt{3}} & \frac{5}{4} & -\frac{5 \sqrt{3}}{4} \\
 0 & -\frac{1}{3} & \frac{\sqrt{3}}{8} & -\frac{1}{8} & \frac{1}{2 \sqrt{3}} & \frac{1}{6} & 0 & 1 &-c_1 & -\frac{35}{24} & -\frac{\sqrt{3}}{2} & -\frac{1}{2} & -\frac{5 \sqrt{3}}{4} & \frac{5}{4} \\
 -\frac{1}{3} & 0 & -\frac{1}{8} & -\frac{\sqrt{3}}{8} & \frac{1}{6} & -\frac{1}{2 \sqrt{3}} & 1 & 0 & -\frac{35}{24} & c_1 & -\frac{1}{2} & \frac{\sqrt{3}}{2} & \frac{5}{4} & \frac{5 \sqrt{3}}{4} \\
 \frac{\sqrt{3}}{2} & -\frac{1}{2} & \frac{1}{2 \sqrt{3}} & \frac{1}{6} & 0 & \frac{1}{4} & -\frac{1}{2 \sqrt{3}} & \frac{1}{6} & -\frac{\sqrt{3}}{2} & -\frac{1}{2} & 0 & \frac{35}{12} & 0 & -\frac{5}{2} \\
 -\frac{1}{2} & -\frac{\sqrt{3}}{2} & \frac{1}{6} & -\frac{1}{2 \sqrt{3}} & \frac{1}{4} & 0 & \frac{1}{6} & \frac{1}{2 \sqrt{3}} & -\frac{1}{2} & \frac{\sqrt{3}}{2} & \frac{35}{12} & 0 & -\frac{5}{2} & 0 \\
 -\frac{\sqrt{3}}{2} & -\frac{1}{2} & \frac{\sqrt{3}}{2} & -\frac{1}{2} & 0 & 1 & -\frac{\sqrt{3}}{2} & -\frac{1}{2} & \frac{\sqrt{3}}{2} & -\frac{1}{2} & 0 & 1 & 0 & 0 \\
 -\frac{1}{2} & \frac{\sqrt{3}}{2} & -\frac{1}{2} & -\frac{\sqrt{3}}{2} & 1 & 0 & -\frac{1}{2} & \frac{\sqrt{3}}{2} & -\frac{1}{2} & -\frac{\sqrt{3}}{2} & 1 & 0 & 0 & 0 \\
\end{array}
\right)}
\end{multline*}

with $c_1 =  \frac{1}{2 \sqrt{3}}-\frac{13 \sqrt{3}}{8}$.

Eigenvalues are $0$ (with multiplicity 4), $\pm i \frac{\sqrt{35}}{4\pi}$ (each with multiplicity 2), $\pm \frac{2}{\pi}$ (each with multiplicity 2) and $\pm \frac{9}{4\pi}$, so on can let $\lambda_0 = \frac{9}{4\pi}$. Computations show that $$\kappa_2 =  \sqrt{\rho\big( \Diff f(Z^*) [\Diff f(Z^*)]^t\big)}+o(1) = \frac{5\sqrt{7}}{2}+o(1).$$
Finally, it is easy to check that it is possible to choose $\nu$ such that relation \eqref{eq:nu} holds, for $\eps$ small enough as soon as
\begin{equation*}
    \nu > \frac{12 + 5\sqrt{7}}{9}+1.
\end{equation*}
Observing that $\frac{12 + 5\sqrt{7}}{9}+1 < 4$, we can conclude that we can construct $\omega_0$ satisfying \eqref{hyp:omega0} with $N=7$ and $\nu =4$ for $\alpha =1$ such that 
\begin{equation*}
    \tau_{\eps,\beta} \le \xi_1 |\ln\eps|,
\end{equation*}
for any $\xi_1 > 4\pi\frac{1-\beta}{9}$.

In order to conclude the proof of Proposition~\ref{prop:N=7}, we observe that $a_N$ and thus $\frac{\kappa_1}{C_\alpha}$, $\frac{1}{C_\alpha}\Diff f(Z^*)$ and thus $\frac{\lambda_0}{C_\alpha}$ and $\frac{\kappa_2}{C_\alpha}$ are all depending continuously on $\alpha$. Therefore, it holds true for $\eps$ small enough that
\begin{equation*}
   \frac{2}{5-\alpha} \left( (1+\alpha)\frac{\kappa_1+\kappa_2}{\lambda_0} + 4-2\alpha\right) < 4.
\end{equation*}
at least on a small interval $[1,\alpha_0)$. This ends the proof.

\vspace{3mm}

\begin{figure}
    \centering
    \includegraphics[width=0.45\textwidth]{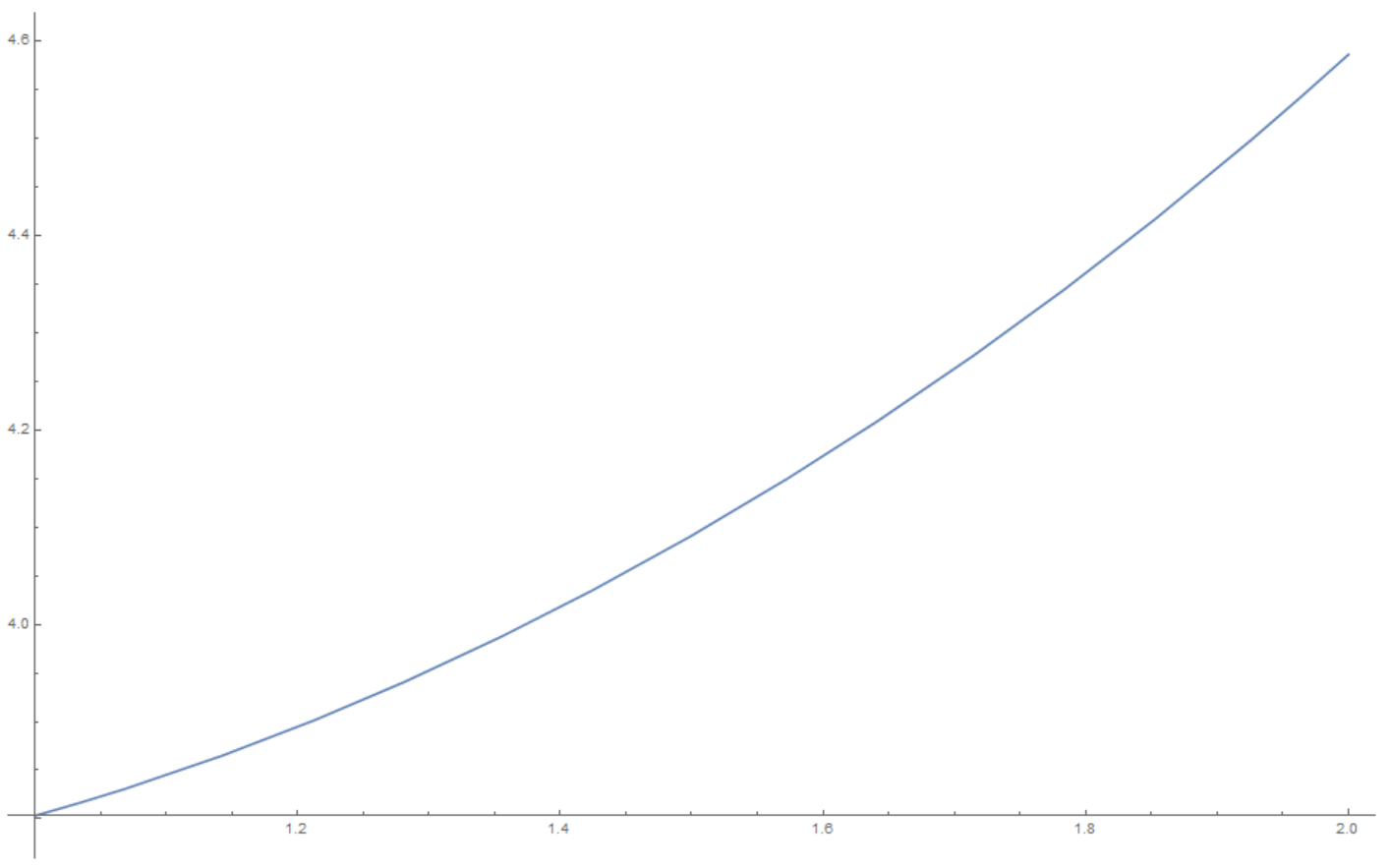} \hspace{5mm}
    \includegraphics[width=0.45\textwidth]{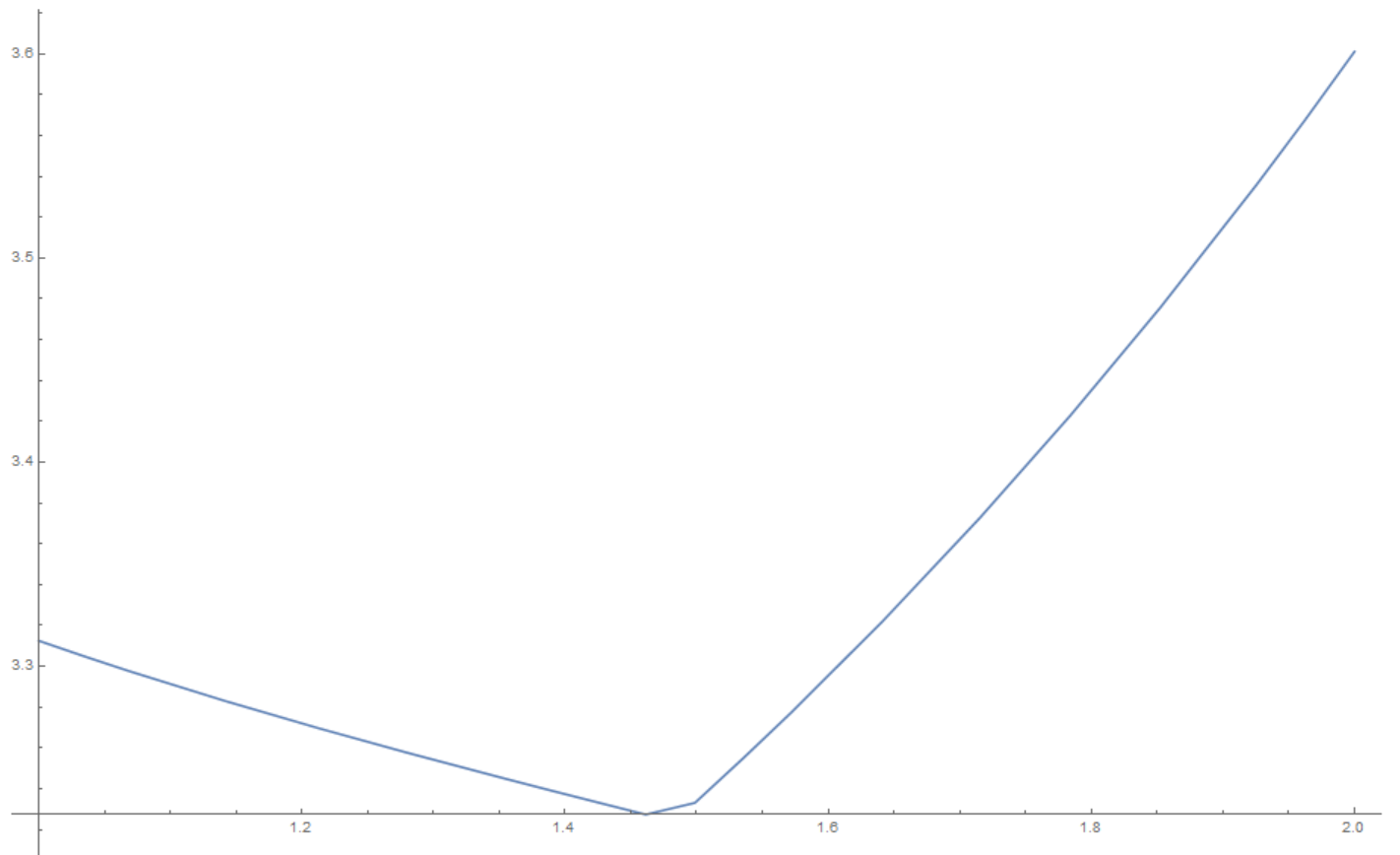}
    \caption{Plot of $h(\alpha)$ in the range $\alpha \in [1,2)$, with $N=7$ (left) and $N=9$ (right).}
    \label{fig:plot_nu_N=7}
\end{figure}

For a general value of $\alpha$, the coefficients of the matrix $\Diff f(Z^*)$ are too complicated to compute mathematically the eigenvalues. However, we use Wolfram Mathematica~\cite{Mathematica} to plot the map
\begin{equation*}
    h: \alpha \mapsto  \frac{2}{5-\alpha}\left( (1+\alpha)\frac{\ds\max_{1\le i \le N}\sum_{j\neq i}\frac{|a_j|}{|z_i^*-z_j^*|^{\alpha+1}}  + \sqrt{\max \big|\mathrm{Eigenvalues}\big[\Diff f(Z^*)[\Diff f(Z^*)]^t\big] \big|}}{\max\big(\Re(\mathrm{Eigenvalues}[\Diff f(Z^*)]\big)} + 4-2\alpha\right)
\end{equation*}
to obtain Figure~\ref{fig:plot_nu_N=7}, which shows that letting $N=9$, we have that $\frac{2}{5-\alpha} \left( (1+\alpha)\frac{\kappa_1+\kappa_2}{\lambda_0} + 4-2\alpha\right)<4$ for every $\alpha \in [1,2)$. Therefore, we have very strong numerical evidence of the fact that a construction is possible for every $\alpha \in [1,2)$ with $\nu < 4$ (thus in particular with $\nu =4$). Please keep in mind that our method does not yield optimal constants. In particular, $\nu =4$ is not optimal.

\paragraph{Acknowledgments.}
\text{ }

The author whishes to acknowledge useful discussions with Thierry Gallay, Pierre-Damien Thizy and Mickaël Nahon. This work was partly conducted when the author was working at the Université Claude Bernard Lyon 1, Institut Camille Jordan.

\bibliographystyle{plain}
\bibliography{base2}

\end{document}